\documentclass[12pt,a4paper]{article}

\newcommand{\bibdir}{../bibtex}

\usepackage{latexsym}
\usepackage{amsfonts}
\usepackage{amssymb}
\usepackage{ifthen}
\usepackage{theorem}
\usepackage{enumerate}
\usepackage{graphics}
\usepackage{color}

\theoremstyle{plain}

\newtheorem{theorem}{Theorem}[section]
\newtheorem{theorem*}[blank]{Theorem}

\newtheorem{conjecture*}[blank]{Conjecture}
\newtheorem{corollary}[theorem]{Corollary}
\newtheorem{corollary*}[blank]{Corollary}
\newtheorem{lemma}[theorem]{Lemma}
\newtheorem{lemma*}[blank]{Lemma}

\newtheorem{proposition*}[blank]{Proposition}

{\theorembodyfont{\rmfamily} 
  \newtheorem{definition}[theorem]{Definition}
  \newtheorem{definition*}[blank]{Definition}
  
  \newtheorem{category*}[blank]{Category}
  
  \newtheorem{functor*}[blank]{Functor}
  
  \newtheorem{algorithm*}[blank]{Algorithm}
  \newtheorem{remark}[theorem]{Remark}
  \newtheorem{remark*}[blank]{Remark}
  
  \newtheorem{notation*}[blank]{Notation}
  
  \newtheorem{problem*}[blank]{Problem}
  
  \newtheorem{question*}[blank]{Question}
  
  \newtheorem{terminology*}[blank]{Terminology}

    \newtheorem{example}[theorem]{Example}
    \newtheorem{example*}[blank]{Example}
}

\newcommand{\isz}{\setlength{\itemsep}{0pt}}
\newcommand{\ifnull}[3]{\def\nullstring{}\def\teststring{#1}\ifx\nullstring\teststring#2\else#3\fi}
\newenvironment{proof}[1][]{\begin{trivlist}\item\ifnull{#1}{\noindent{\sc Proof:\ }}{\noindent{\sc Proof of #1:\ }}}{\hspace*{\fill}$\Box$\vspace*{\baselineskip}\end{trivlist}}
\newcounter{enum}

\newcommand{\figref}[1]{Figure~\ref{fig:#1}}

\newcommand{\figrefpart}[2]{Figure~\ref{fig:#1}(#2)}

\newcommand{\figscale}{1}

\newcommand{\fig}[2]{
  \begin{figure}[htb]\centering
    \ifnull{#1}{}{\scalebox{\figscale}{\includegraphics{#1.eps}}}
    \ifnull{#2}{}{\caption{#2}}
    \label{fig:#1}
  \end{figure}
}

\newcommand{\R}{\mathbb{R}}

\newcommand{\Z}{\mathbb{Z}}

\renewcommand{\th}{\ensuremath{^\mathrm{th}}}

\renewcommand{\emptyset}{\varnothing}
\renewcommand{\geq}{\geqslant}
\renewcommand{\leq}{\leqslant}

\newcommand{\card}[1]{\#(#1)}
\newcommand{\dirlim}[1]{\ifthenelse{\equal{#1}{}}{\underrightarrow{\lim}}{\underrightarrow{#1}}}

\newcommand{\fto}{\longrightarrow}
\newcommand{\exto}{\fto}

\newcommand{\homotopic}{\sim}

\newcommand{\hniel}{h_\mathit{niel}}
\newcommand{\htop}{h_\mathit{top}}

\newcommand{\p}{^\prime}

\newcommand{\tendsto}{\rightarrow}

\newcommand{\pre}[1]{\textrm{Pre-}#1}

\newcommand{\us}{{U/S}}

\newcommand{\edge}[1]{E(#1)}
\newcommand{\control}[1]{C(#1)}
\newcommand{\vertex}[1]{V(#1)}

\newcommand{\init}[1]{\imath(#1)}
\newcommand{\final}[1]{\imath(\bar #1)}

\renewcommand{\card}[1]{|#1|}
\newcommand{\cutf}[1]{\mathcal{C}{#1}}
\newcommand{\cutm}[1]{\mathcal{C}_{{#1}^U}{M}}
\newcommand{\cutt}[1]{\mathcal{C}{#1}}
\newcommand{\cutp}[1]{(\mathcal{C}_{{#1}^U}{M},\mathcal{C}_{{#1}^U}{{#1}^S})}

\newcommand{\normalmargins}{
  \setlength{\oddsidemargin}{-0.4mm}
  \setlength{\evensidemargin}{-0.4mm}
  \setlength{\textwidth}{160mm}

  \setlength{\topmargin}{-0.4mm}
  \setlength{\headheight}{4mm}
  \setlength{\headsep}{7mm}
  \setlength{\textheight}{225mm}
  \setlength{\footskip}{11mm}

  \addtolength{\topsep}{-\parskip}
  \addtolength{\partopsep}{-\parskip}

  \setlength{\headheight}{0mm}
  \setlength{\headsep}{0mm}
  \setlength{\textheight}{236mm}
}

\renewcommand{\figscale}{1}
\normalmargins

\begin{document}
\sloppy
%\doublespace

%****************************************************************************************************************************************************************************

\title{Entropy-minimising models of surface diffeomorphisms relative to homoclinic and heteroclinic orbits}
\author{Pieter Collins
  \thanks{This work was partially funded by Leverhulme Special Research Fellowship SRF/4/9900172.}
  \\Centrum voor Wiskunde en Informatia\\P.O. Box 94079\\1090 GB Amsterdam\\The Netherlands\\pieter.collins@cwi.nl
}
%\date{Draft \today}
\date{\today}
\maketitle

\begin{abstract}
In the theory of surface diffeomorphisms relative to homoclinic and heteroclinic orbits, it is possible to compute a one-dimensional representative map for any irreducible isotopy class.
The topological entropy of this graph representative is equal to the growth rate of the number of essential Nielsen classes of a given period, and hence is a lower bound for the topological entropy of the diffeomorphism.
In this paper, we show that this entropy bound is the infemum of the topological entropies of diffeomorphisms in the isotopy class, and give necessary and sufficient conditions for the infemal entropy to be a minimum.
\end{abstract}


\newpage

%{\setlength{\parskip}{0pt} \tableofcontents}

\newpage

%****************************************************************************************************************************************************************************

\section{Introduction}
\label{sec:introductions}

In the Nielsen-Thurston theory of surface diffeomorphisms \cite{CassonBleiler88}, it is possible to find a diffeomorphism in each isotopy class which minimises both the number of periodic orbits of each period, and the topological entropy \cite{Handel85ERGTD,Fathi1990ERGTD,Jiang1996,Boyland1999AMS}.
A constructive proof can be given \cite{BestvinaHandel1995TOPOL,FranksMisiurewicz1993}; a one-dimensional representative of the diffeomorphism called a \emph{train
  track} is computed and the Thurston-minimal model is easily constructed from this.

The case of diffeomorphisms relative to homoclinic orbits to a saddle point is more complicated.
It is well-known that any diffeomorphism with a transverse homoclinic point
has positive topological entropy, but it is easy to construct examples of such
diffeomorphisms with arbitrarily small entropy.
Therefore, it may not even be possible to find a diffeomorphism realising the lower bound of the entropy.

In this paper, we discuss the relation between the \emph{Nielsen entropy} \cite{Jiang83}, which measures the growth rate of the number of essential Nielsen classes of periodic points, and the infemal entropy in the isotopy class.
The Nielsen entropy is equal to the topological entropy of the \emph{graph representative}, which can be computed using the algorithms described in \cite{Collins2002INTJBC,CollinsPPDYNSYS}.
The key step is to represent the isotopy class relative to the homoclinic/heteroclinic orbits using a trellis, which is a subset of the homoclinic tangle of the saddle periodic points.
We show that under mild non-degeneracy conditions, the Nielsen entropy is equal to the infemum of the topological entropies in the isotopy class.
This result can be considered an optimality result for the computational techniques of the trellis theory.

The proof that the Nielsen entropy is infemal involves the construction of uniformly-hyperbolic diffeomorphisms realising the entropy bound arbitrarily closely.
In certain cases, it is possible to construct a uniformly-hyperbolic diffeomorphism realising the entropy bound.
We show that the only obstruction to the construction of a uniformly hyperbolic diffeomorphisms realising the entropy bound is the existence of \emph{almost wandering segments}; intervals of (un)stable manifold which must contain intersection points under some iterate of any representative diffeomorphism, but for which no bound on the number of iterates exists. 
Hence the Nielsen entropy is a minimum if and only if there are no almost wandering segments.
Under these conditions, we also show that if the system \emph{irreducible}, then any isotopy which removes any of the intersection points of the trellis also reduces the Nielsen entropy.

The conditions on the existence of entropy minimisers is delicate, and relies on the class of system under consideration.
We only consider the problem in the case of homoclinic orbits to a saddle point in which the local stable and unstable branches are known, such as a hyperbolic saddle periodic orbit of a diffeomorphism which are the typical case in applications.
If we instead consider systems for which we allow periodic points modelled on \emph{$n$-prong singularities}, then it is immediate that there are no almost wandering segments, and it is always possible to introduce extra branches ending in attracting or repelling periodic orbits and hence construct a homeomorphism realising the entropy bound.
A similar situation occurs if we only require that the stable and unstable curves of the trellis are, respectively, positively and negatively invariant.

\pagebreak

%****************************************************************************************************************************************************************************

\section{Trellises, curves and graphs}
\label{sec:tralliscurvegraph}

We introduce the notion of trellis and trellis type, and relate these to the well-known notion of homoclinic/heteroclinic orbit and tangle.
We give a number of definitions which allow us to describe some of the basic properties of trellises.
We then give some definitions concerning homotopy and isotopy classes of curves.
Finally, we review the definitions of (thick) graphs, and show how to associate a graph representative to each irreducible trellis type.

%----------------------------------------------------------------------------------------------------------------------------------------------------------------------------

\subsection{Tangles and Trellises}
\label{sec:tangletrellis}

A \emph{tangle} is the figure formed by the stable and unstable manifolds of a collection of periodic saddle orbits $P$ of a surface diffeomorphism $f$.
We denote the stable and unstable sets of $P$ by $W^S(f;P)$ and $W^U(f;P)$ respectively.
If $P$ consists of a single point, then $W$ is a \emph{homoclinic} tangle, otherwise $W$ is a \emph{heteroclinic} tangle.
\begin{definition}[Trellis]
Let $f$ be a diffeomorphism of a surface $M$ with a finite invariant set $P$ of hyperbolic saddle points.
A \emph{trellis} for $f$ is a pair $T=(T^U,T^S)$, where $T^U$ and $T^S$ be subsets of $W^U(f;P)$ and $W^S(f;P)$ respectively such that: 
\begin{enumerate}
\item $T^U$ and $T^S$ both consist of finitely many compact intervals with non-empty interiors,
\item $f(T^U)\supset T^U$ and $f(T^S)\subset T^S$.
\end{enumerate}
We denote the set of periodic points of $T$ by $T^P$, and the set of intersections of $T^U$ and $T^S$ by $T^V$.
\end{definition}
We write $(f;T)$ to denote the pair consisting of a diffeomorphism $f$ and a trellis $T$ for $f$.
We use the notation $T^{\us}$ in statements which hold for both stable and unstable manifolds.
If the diffeomorphism $f$ under consideration is clear from the context, we abbreviate $T^{\us}(f;P)$ to $T^{\us}(P)$,
 and if the periodic point set $P$ is also clear, we simply write $T^{\us}$.
If $p$ is a point of $P$, then the stable and unstable curves passing through $p$ are denoted $T^S(f;p)$ and $T^U(f;p)$ respectively.
We denote the closed interval in $T^{\us}$ with endpoints $a$ and $b$ by $T^{\us}[a,b]$ and the closed interval by $T^{\us}(a,b)$.

An \emph{intersection point} of a tangle $T$ is a a point in $T^U\cap T^S$, and is a point of a homoclinic or heteroclinic orbit.
An intersection point $q\in T^U(p_1)\cap T^U(p_2)$ is a \emph{primary intersection point} or \emph{pip} if $T^U(p_1,q)$ and $T^S(q,p_2)$ are disjoint.
An intersection point $q$ is \emph{transverse} or \emph{tangential} according to whether $T^U$ and $T^S$ cross transversely or tangentially at $q$
An \emph{end intersection} of a trellis is an intersection which does not lie in two segments of $T^\us$ (so lies in at least one end of $T^\us$).

An \emph{arc} of $T^\us$ is a closed subinterval of $T^\us$ with endpoints in $T^V$.
A \emph{segment} of $T^\us$ is an arc of $T^\us$ with no topologically transverse intersection points in its interior.
An \emph{end interval} of a trellis is a subintervals of $T^\us$ which does not lie in any segment.

A \emph{branch} of a trellis $T$ is a set of the form $ \{p\}\cup W^{\us}(p,\pm\infty)\cap T^{\us} . $
A branch is \emph{trivial} if it does not intersect any other branch.
Note that points $q_1$ and $q_2$ lie in the same branch of $T^\us$ if and only if $f^{\mp n}(q_1)$ and $f^{\mp n}(q_2)$ lie in the same component of $T^{\us}\setminus P$ for some $n$.
Note that a branch need not be connected.
 
An \emph{open region} of $T$ is a component of $M\setminus(T^U\cup T^S)$. 
A \emph{(closed) region} is the closure of an open region, and hence includes the stable and unstable boundary segments.
A \emph{(closed) domain} of $T$ is an closed subset of $M$ with boundary in $T^U\cup T^S$;
 note that we do not require domains to be simply-connected, though this is the case we will usually consider.
We will denote domains by the letter $D$, and regions by the letter $R$.
Note that we do not require domains or regions to be simply-connected.
There are two special types of region which play an important role, namely \emph{bigons} and \emph{rectangles}.
A \emph{bigon} is a region which is a topological disc bounded by one stable and one unstable segment with internal angles less than $\pi$.
A \emph{rectangle} is a region which is a topological disc bounded by two stable and two unstable segments with internal angles less than $\pi$.
Bigons play a similar role in the trellis theory as do critical points in the kneading theory, and may be called \emph{critical} regions; rectangles can be foliated by stable and unstable leaves, and may be called \emph{regular} regions.

In general, we do not have much control over the geometry of a trellis.
The exception is near a point of $T^P$, where the dynamics are conjugate to a linear map.
It is important to consider how the nontrivial branches $T$ divide the surface in a neighbourhood of $T^P$.
\begin{definition}[Quadrant, secant and coquadrant]
Let $T$ be a trellis, and $\overline{T}$ be the nontrivial branches of $T$, and $p$ be a point of $T^P$.
Then a \emph{sector} of $T$ at $p$ is a local component of a region of $\overline{T}$ in a neighbourhood of $p$.
A sector is a \emph{quadrant} $Q$ if it intersects no trivial branches, in which case the boundary includes a single unstable branch $T^U(Q)$ and a single stable branch $T^S(Q)$ with internal angle less than $\pi$.
A sector is an \emph{secant} if it intersects a single trivial branch; if this is a branch of $T^U$ the secant is \emph{attracting}, and if this is a branch of $T^S$, the secant is repelling.
A sector is a \emph{coquadrant} if it intersects two trivial branches, and so the nontrivial branches $T^U(Q)$ and $T^S(Q)$ subtend an angle greater than $\pi$ in $Q$.
\end{definition}

\fig{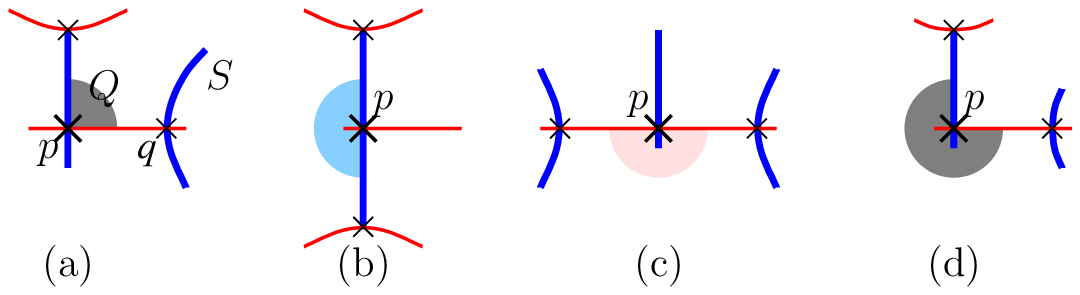}{Sectors at $p$. (a) depicts a quadrant $Q$, (b) an attracting secant, (c) a repelling secant and (d) a coquadrant.}

Quadrants, secants and coquadrants are shown in \figref{quadrantsecant}.
The region containing a quadrant $Q$ is denoted $R(Q)$.
The \emph{image} of a quadrant $Q$ is the quadrant containing $f(Q)$; note that $T^\us(f(Q))$ is the same branch as $f(T^\us(Q))$.
A stable segment $S$ with endpoint $q\in T^U(Q)$ intersects lies on the \emph{$Q$-side} of $T^U$ if locally $S$ lies on the same side of $T^U(Q)$ at $q$ as $T^S(Q)$ does at $p$.
Similarly, an unstable segment $U$ with endpoint $q\in T^S(Q)$ intersects lies on the \emph{$Q$-side} of $T^S$ if locally $U$ lies on the same side of $T^S(Q)$ at $q$ as $T^U(Q)$ does at $p$.
The segment $S$ in \figrefpart{quadrantsecant}{a} lies on the $Q$-side of $q$.

By an \emph{isotopy relative to $T$}, we mean an isotopy $f_t$ such that $T$ is a tangle for $f_t$ for all $t$.
We say $(f_0;T_0)$ and $(f_1;T_1)$ are \emph{conjugate} if there is a homeomorphism $h$ such that $h\circ f_0=f_1\circ h$ and $h(T_0^{\us})=T_1^{\us}$.
These relations allow us to define the equivalence classes of trellis map which will be our primary object of study.
\begin{definition}[Trellis mapping class and trellis type]
The \emph{trellis mapping class} $([f];T)$ is the set of all pairs $(\tilde{f};T)$ for which $\widetilde{f}$ is isotopic to $f$ relative to $T$.
The \emph{trellis type} $[f;T]$ is the set of all pairs $(\tilde{f},\widetilde{T})$ which are conjugate to a map isotopic to $f$ relative to $T$.
We consider $(f_0;T_0)$ and $(f_1;T_1)$ equivalent if $[f_0;T_0]=[f_1;T_1]$.
\end{definition}
A trellis type $[\widehat{f};\widehat{T}]$ is an iterate of the trellis type $[f;T]$ if $(\widehat{f},T)\in[f;T]$ and $\widehat{T}$ is an $\widehat{f}$-iterate of $T$.

A surface diffeomorphism $f$ with a periodic saddle orbit has infinitely many trellises, which are partially ordered by inclusion.
Taking a smaller trellis gives a \emph{subtrellis}, and a larger trellis a \emph{supertrellis}.
Then $T$ is a \emph{subtrellis} of $\widehat{T}$ if $T^{\us}\subset \widehat{T}^{\us}$, and $T$ is a subtrellis of $W$ if $T^{\us}\subset W^{\us}$.
We say $\widehat{T}$ is a \emph{supertrellis} of $T$ and $\widehat{W}$ is a \emph{supertangle} of $T$.
Similarly, the trellis type $[f;T]$ is a subtrellis of $[f;\widehat{T}]$.
Of particular importance are those supertrellises which can be obtained by iterating segments or branches.
A trellis $\widehat{T}$ is an \emph{$f$-iterate} of $f$ if there exist positive integers $n_u$ and $n_s$ such that $\widehat{T}^U=f^{n_u}(T^U)$ and $\widehat{T}^S=f^{-n_u}(T^S)$.
A trellis $\widehat{T}$ is an \emph{$f$-extension} of $T$ if there exists $n$ such that
  \[ T^U\subset \widehat{T}^U\subset f^n(T^U) \quad \textrm{and} \quad T^S\subset \widehat{T}^S\subset f^{-n}(T^S) . \]
An iterate/extension is a \emph{stable iterate}/\emph{extension} if $\widehat{T}^U=T^U$ and an \emph{unstable iterate}/\emph{extension} if $\widehat{T}^S=T^S$.

The main difference between extensions and supertrellises is that $\widehat{T}^P=T^P$ for an extension, but $T^P$ may be a strict subset of $\widehat{T}^P$ for a supertrellis.
This difference makes the analysis of supertrellises slightly more complicated than that of extensions.

If $S$ is a segment of $T^S$ and $f^{-1}(S)$ is not a subset of $T^S$, we can take a stable $f$-extension $\widehat{T}$ with $\widehat{T}^S=T^S\cup f^{-1}(S)$.
We say this extension is formed by taking a backward iterate of $S$.
In a similar way, we can form extensions by taking backward iterates of arcs of $T^S$, or forward iterates of arcs of $T^U$.
For these extensions, the endpoints of any curve of $\widehat{T}$ lie in the set $X=\bigcup_{n=-\infty}^{\infty}f^n(T^V)$.

We say a supertrellis $(\widetilde{f},\widehat{T})$ is a \emph{minimal supertrellis} of $([f];T)$ if $\widehat{T}$ has minimal number of intersections among all supertrellises with endpoints in the same segment. An equivalent formulation of a minimal supertrellis is that every bigon contains a point of the set $X=\bigcup_{n\in\Z}\widetilde{f}^n(T^V)$. 
A minimal supertrellis which is an extension or iterate is, respectively, a \emph{minimal extension} or a \emph{minimal iterate}.

The most important dynamical feature of a trellis type is its \emph{entropy}.
\begin{definition}[Entropy]
The \emph{entropy} of a trellis type $[f;T]$, denoted $\htop[f;T]$ is the infemum of the topological entropies of diffeomorphisms in $[f;T]$; that is
\[ \htop[f;T]=\inf\{\htop(\widehat{f}):\widehat{f}\in[f;T]\} \;. \]
If this infemum is not a minimum, we sometimes write $\htop[f;T]=\inf\{\htop(\widehat{f})\}+\epsilon$.
\end{definition}

\begin{example}[The type-$3$ trellis]
\fig{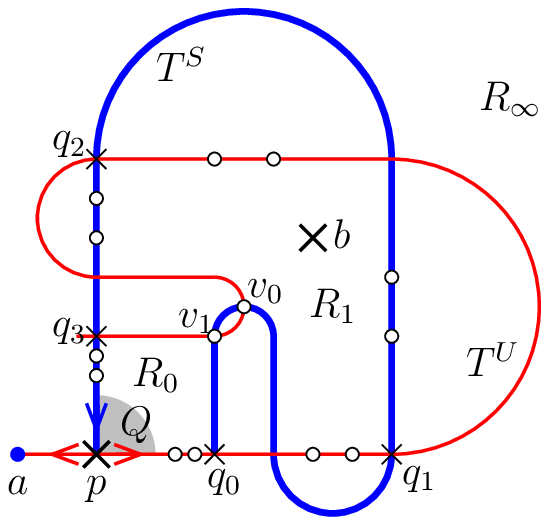}{The trellis type $[f_3;T_3]$}

The type-$3$ trellis is depicted in \figref{type3}.
It is formed by subsets of the stable and unstable manifolds of the direct saddle fixed point (i.e. the saddle point with positive eigenvalues)
 of the H\'enon map for certain parameter values.
The stable and unstable sets are subsets of the stable and unstable manifolds of the saddle fixed point $p$.
The branches of this trellis are connected and the all intersection points are transverse.
The points $q_0$, $q_1$, $q_2$ and $q_3$ are primary intersection points on a single homoclinic orbit.
The orbits of the intersection points $v_0$ and $v_1$ are shown in white dots.
These orbits are called \emph{forcing orbits} for the type-$3$ trellis; 
 it is essentially impossible to remove any intersection points of the trellis by an isotopy without first removing $v_0$ and $v_1$ in a homoclinic bifurcation.
Throughout this paper we follow the convention of showing forcing orbits as which dots, and primary intersection points as crosses.

There are ten regions, an unbounded region $R_\infty$, five bigons, three rectangles and a hexagon.
The quadrant $Q$ is contained in the region $R_0$.
Under the Smale horseshoe map, there is a Cantor set of nonwandering points contained in the (closed) regions $R_0$ and $R_1$, including a fixed point $b$ in $R_1$.
All other points are wandering except for an attracting fixed point at the end of one unstable branch of $T^U(p)$
 in $R_\infty$.
The topological entropy of the Smale horseshoe map is $\log\lambda_{\max}$, where $\lambda_{\max}$ is the largest root of the polynomial $\lambda^3-\lambda^2-2$.
Any diffeomorphism $f$ with this trellis type must have topological entropy $\htop(f)\geq\log\lambda_{\max}\approx0.528$.

Notice that one of the unstable branches of the trellis ends in an attracting fixed point,
 and one of the stable branches is non-existent.
Both of these are therefore trivial branches.
\end{example}

Given a diffeomorphism $f$ with a trellis $T$, we can obtain a canonical map of pairs by \emph{cutting} along $T^U$, as described in \cite{Collins1999AMS}.
The topological pair obtained by cutting along the unstable curve is denoted $\cutt{T}=\cutp{T}$.
The diffeomorphism $f$ lifts for a map $\cutf{f}$ on $\cutt{T}$.
Notice that the pair $\cutt{T}=\cutp{T}$ contains the pair $(M\setminus T^U,T^S\setminus T^U)$ as an open subset which is invariant under $\cutf{f}$.
Indeed, $\cutm{T}$ can be regarded as a natural compactification of $M\setminus T^U$, and the homotopy properties of $\cutp{T}$ and $(M\setminus T^U,T^S\setminus T^U)$ are essentially the same.
In order that the stable set be well-behaved after cutting, we only consider \emph{proper trellises}, for which $\partial T^S\subset T^U$ and $\partial T^U\cap T^S=\emptyset$.
In other words, a trellis is proper if the endpoints of intervals in $T^S$ lie in $T^U$, but the endpoints of intervals in $T^U$ do not lie in $T^S$.

\begin{definition}[Well-formed trellises]
We say a proper trellis $T=(T^U,T^S)$ for a diffeomorphism $f$ is \emph{well-formed} if every component of $T^U\cup f(T^S)$ contains a point of $T^P$.
\end{definition}
Restricting to well-formed trellises will be important when considering the graph representative defined in Section~\ref{sec:graphmaps}, as it gives a necessary condition for the topological entropy of the graph representative to be an optimal entropy bound for the trellis mapping class.

%----------------------------------------------------------------------------------------------------------------------------------------------------------------------------

\subsection{Attracting and repelling regions}
\label{sec:attractor}

As well as saddle periodic orbits, codimension-zero sets of diffeomorphisms also have attracting and repelling periodic orbits.
Clearly, a stable segment cannot be in the basin of attraction of an attracting orbit, and neither can an unstable segment be in the basin of a repelling orbit.
However, it is possible for the interior of an entire region to lie in the basin of an attracting or repelling periodic orbit.
Such regions are called \emph{stable} or \emph{unstable} regions.
\begin{definition}[Stable and unstable regions]
Let $([f];T)$ be a trellis mapping class. A domain $D$ of $T$ is \emph{stable} if there is a diffeomorphism $\widetilde{f}$ in $([f];T)$ and a periodic orbit $P$ of $f$ such that
 $\widetilde{f}^n(x)\tendsto P$ as $n\tendsto\infty$ for all $x$ in the interior of $D$.
If $D$ contains a point of $P$, then $D$ is \emph{attracting}.
Similarly, a domain $D$ is \emph{unstable} if there is a diffeomorphism $\widetilde{f}$ in $([f];T)$ and a periodic orbit $P$ of $f$
 such that $\widetilde{f}^{-n}(x)\tendsto P$ as $n\tendsto\infty$ for all $x$ in the interior of $D$, and if $D$ contains a point of $P$, then $D$ is \emph{repelling}.
\end{definition}

\fig{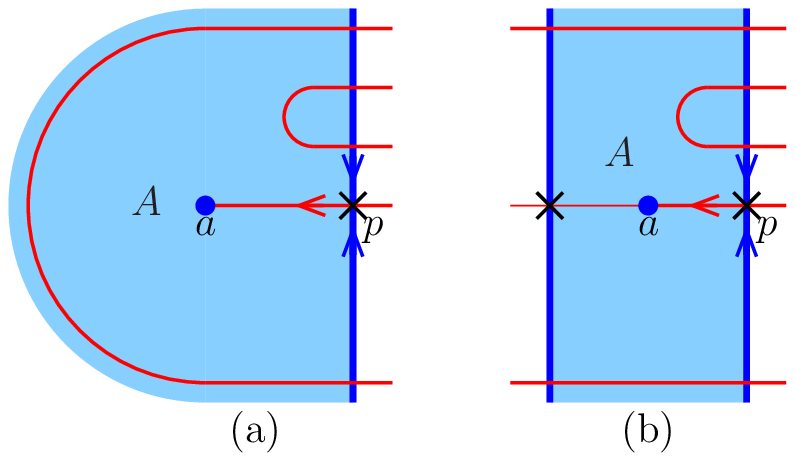}{Attracting domains.}
Two examples of attracting domains are shown in \figref{attractor}.
All points in the shaded domains are in the basin of attraction of the stable periodic point.
Note that $T^S$ is disjoint from the interior of $A$.
Repelling domains are similar.
Stable and unstable domains are wandering sets, unless they contain the attracting or repelling periodic orbit.

\fig{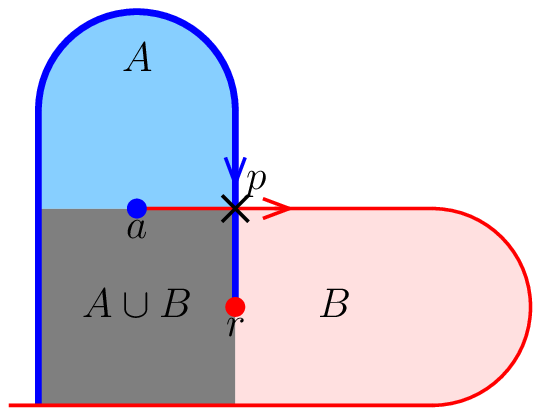}{An unbounded region.}
A region which contains an attractor and a repellor is called \emph{unbounded} region, by analogy with the unbounded region of the planar Smale horseshoe, for which the repellor is the point at infinity.
\figref{infinity} shows an unbounded region $R$ near a fixed point $p$. 
All points in $A$ lie in the basin of attraction of an attracting fixed point $a$,
 and all points in $B$ lie in the basin of a repelling fixed point $b$.
Since $A$ and $B$ cover the interior of $R$, all points in the interior of $R$ are wandering, apart from $a$ and $b$.

\begin{definition}[Chaotic region]
Let $([f];T)$ be a trellis mapping class.
A region $R$ is \emph{chaotic} if for every diffeomorphism $\widehat{f}$ in $([f];T)$, there exists an integer $n$ such that $f^n(T^U)\cap f^{-n}(T^S)$ contains a point in the interior of $R$.
\end{definition}
We shall see that a chaotic region supports entropy, in the sense that for any diffeomorphism $\widehat{f}\in([f];T)$, there is an ergodic measure $\mu$ with positive entropy such that $\mu(R)>0$.

An open segment which is in the basin of a repelling periodic orbit is also a wandering set, as is an open segment which is in the basin of an attracting periodic orbit.
There are some nonwandering segments such that any finite iterate need not contain an intersection point.
We call such segments \emph{almost wandering}.
\begin{definition}[Almost wandering segment]
\label{defn:almostwandering}
Let $([f];T)$ be a trellis mapping class.
An open segment is \emph{almost wandering} if it is nonwandering, but for any integer $n$ there exists a diffeomorphism $\tilde{f}\in([f];T)$ the $n\th$ iterate of the segment by $\tilde{f}$ contains no intersection points.
\end{definition}
Hence an open stable segment $S$ is almost wandering if for any positive integer $n$, there exists $\widetilde{f}\in([f];T)$ such that $\widetilde{f}^{-n}(S)\cap T^U=\emptyset$, but for any diffeomorphism $\widehat{f}\in([f];T)$, there exists $n$ such that $\widehat{f}^{-n}(S)\cap T^U\neq\emptyset$.
A similar statement holds for unstable segments.

%----------------------------------------------------------------------------------------------------------------------------------------------------------------------------

\subsection{Curves and homotopies}
\label{sec:curves}

Our main tool for studying the geometry, topology and dynamics associated with trellis maps will be to consider curves embedded in pairs $(X,Y)$, which we will usually take to be either the cut surface $\cutp{T}$ of the graph $(G,W)$.
Since the topological pair $\cutt{T}$ is obvious from the trellis, we will usually draw curves in $\cutm{T}$ as curves embedded in the original surface $M$, and,
 wherever possible, ensure these curves are disjoint from $T^U$.

Our curves will be maps in this category of topological pairs, $\alpha:(I,J)\fto(X,Y)$.
The \emph{path} of such a curve $\alpha$ is the set $\alpha(I)$.
For the most part, we are only interested in curves up to homotopy or isotopy, and we always take homotopies and isotopies of curves through maps of pairs.
For Nielsen theory we will always keep the endpoints fixed during the homotopy.
For most other purposes, we only consider curves for which $J\supset\{0,1\}$, so the endpoints lie in $Y$, but may move.
If $\alpha:(I,J)\fto(X,Y)$ is a curve and $J$ contains $\{0,1\}$ we say $\alpha$ has endpoints in $Y$; if $J$ equals $\{0,1\}$ we say $\alpha$ has endpoints only in $Y$.

Reparameterising a curve does not change its path, but may change the set $J$ which maps into $T^S$.
This means that different parameterisations of the same path may not even be comparable under homotopy.
However, we consider different parameterisations of the same curve as equivalent.
\begin{definition}[Equivalence of curves]
Curves $\alpha_1:(I_1,J_1)\fto(X,Y)$ and $\alpha_2:(I_2,J_2)\fto(X,Y)$ are \emph{equivalent} if there is a homeomorphism $h:(I_1,J_1)\fto(I_2,J_2)$ with $h(J_1)=J_2$ such that $\alpha_1\homotopic\alpha_2\circ h$ as curves $(I_1,J_1)\fto(X,Y)$.
The homotopy may be taken relative to endpoints, as appropriate.
\end{definition}
We can also define a partial order on curves in a similar way:
\begin{definition}[Tightening curves]
\label{defn:curvetightening}
Let $\alpha_1:(I_1,J_1)\exto(X,Y)$ and $\alpha_2:(I_2,J_2)\exto(X,Y)$ be exact curves.
We say $\alpha_2$ \emph{tightens} to $\alpha_1$ if there is an injective map $h:(I_1,J_1)\fto(I_2,J_2)$ such that $\alpha_1\homotopic\alpha_2\circ h$ as curves $(I_1,J_1)\fto (X,Y)$.
\end{definition}
If $\alpha_2$ tightens to $\alpha_1$, then $\alpha_2$ can be thought of as more complicated than $\alpha_1$.
It is clear that the tightening relation is a partial order on homotopy classes of curves.
We are especially interested in iterates with minimal number of components intersecting $Y$.
\begin{definition}[Minimal iterate]
Let and $\alpha:(I,J)\exto(X,Y)$ be a simple curve with endpoints in $Y$.
Then a \emph{minimal iterate} of $\alpha$ under $f$ is a curve $\beta$ which is homotopic to $f\circ\alpha$ relative to $J$ and which minimises the number of components of intersections with $Y$.
We let $J\p=\beta^{-1}(Y)$, and consider $\beta$ as a curve $(I,J\p)\fto(X,Y)$.
If we further require that in intersection is isolated whenever possible, then the curve $\beta$ is well-defined up to equivalence, so we obtain a well-defined map $f_{\min}$ on equivalence classes of curves given by $f_{\min}[\alpha]=[\beta]$.
\end{definition}
Note that by $f_{\min}^n[\alpha]$ we mean ${(f_{\min})}^n[\alpha]$, and not ${(f^n)}_{\min}[\alpha]$, which typically has fewer intersections.

%----------------------------------------------------------------------------------------------------------------------------------------------------------------------------

\subsection{Graph maps}
\label{sec:graphmaps}

Our main tool for computing and describing the dynamics forced by a trellis mapping class or trellis type is to relate the trellis map to a graph map.
We will mostly use the same terminology and notation for graphs as in \cite{BestvinaHandel1995TOPOL}, though our terminology for thick graphs is closer to that of \cite{FranksMisiurewicz1993}.

In particular, a \emph{graph} $G$ is a one-dimensional CW-complex with vertices $\vertex{G}$ and edges $\edge{G}$.
The reverse of an edge $e$ is denoted $\bar{e}$.
An \emph{edge-path} is a list $e_1\ldots e_n$ of oriented edges of $G$ such that $\final{e_i}=\init{e_{i+1}}$ for $1\leq i < n$, and an \emph{edge loop} is a cyclically-ordered list of edges.
The \emph{trivial} edge-path contains no edges and is denoted $\cdot$.
An edge-path $e_1\ldots e_n$ \emph{back-tracks} if $e_{i+1}=\bar{e}_i$ for some $i$, otherwise it is \emph{tight}.
A \emph{graph map} $g$ is a self-map of $G$ taking a vertex to a vertex, and edge $e$ to an edge-path $e_1\ldots e_k$
 such that $\init{e_1}=g(\init{e})$ for all directed edges $e$.
The \emph{derivative} map $\partial g$ takes oriented edges to oriented edges or $\cdot$, with $\partial g(e_i)=e_j$
 if $g(e_i)=e_j\ldots$ and $\partial g(e_i)=\cdot$ if $g(e_i)=\cdot$.

We will always consider a graph embedded in a surface by an embedding $i$.
This induces a natural cyclic order $\lhd$ on the oriented edges starting at each vertex.
A pair of edges $(e_1,e_2)$ is a \emph{turn} in $G$ at at vertex $v$ if $v=\init{e_1}=\init{e_2}$ and $e_1\lhd e_2$, so $e_2$ immediately follows $e_1$ in the cyclic order at $v$.
An edge-loop $\pi=\ldots,p_1,p_2,\ldots,p_n,\ldots$ is \emph{peripheral} in $G$ if $(p_{i+1},\bar{p}_i)$ is a turn in $G$ for all $i$.

A graph embedded in a surface can be ``thickened'' to obtain a \emph{fibred surface} or a \emph{thick graph} with \emph{thick vertices} and \emph{thick edges}, the latter fibred by stable and unstable leaves.
The inverse of this thickening is achieved by collapsing thick vertices to points, and thick edges to edges.
A \emph{thick graph map} is a injective map of a thick graph which is a strict contraction on the set of thick vertices, and maps thick edges so as to preserve the stable leaves.
A graph map is \emph{embeddable} if it can be obtained from a thick graph map by collapsing the stable leaves.
The \emph{peripheral subgraph} $P$ of $g$ is a maximal invariant subset of $G$ consisting of simple peripheral loops.
Edges of $P$ are called \emph{peripheral edges}.
If $g^n(e)\subset P$ for some $n$, then $e$ is \emph{pre-peripheral}.
The set of pre-peripheral edges is denoted $\pre{P}$, and contains $P$.

The \emph{transition matrix} of a graph map $g$ is the matrix $A=(a_{ij})$ where $a_{ij}$ counts the number of times the edge $e_j$ appears in the image path of edge $e_i$.
The largest eigenvalue of $A$ is the \emph{growth rate} $\lambda$ of $g$, and the logarithm of the growth rate gives the topological entropy of $g$.
A \emph{length function} on $G$ is a strictly-positive function $l:\edge{G}\fto\R$.
The \emph{length} of an edge-path $e_1e_2\ldots e_n$ is defined to be $l(e_1e_2\ldots e_n)=\sum_{i=1}^nl(e_i)$.
If $g$ is a graph map with topological entropy $\lambda$, then for any $\epsilon>0$ there is a length function with $l(g(e))<(\lambda+\epsilon)l(e)$ for any edge $e$. 

We relate a trellis mapping class to a graph map via the map $\cutf{f}$ of the topological pair $\cutt{T}$ obtained by cutting along the unstable curves.
The result is a map of a topological pair $(G,W)$ where $W$ is a finite subset of $G$.
\begin{definition}[Controlled graph map]
Let $(G,)$ be a topological pair where $G$ is a graph and $W$ is a finite subset of $G$.
The edges of $G$ containing points of $W$ are called \emph{control edges}.
A graph map $g:(G,W)\fto(G,W)$ is a \emph{controlled graph map} if $g(z)$ is a control edge whenever $z$ is a control edge.
\end{definition}

A vertex of $G$ which is the endpoint of a control edge is called a \emph{control vertex}.
All other vertices are called \emph{free vertices}, and edges which are not control edges are called \emph{free edges}.
Free edges which are neither peripheral nor pre-peripheral are called \emph{expanding edges.}

A controlled graph is \emph{proper} if it is connected and each free vertex has valence at least $3$.
A proper controlled graph has at most $3\card{\control{G}}-3\chi(G)$ edges, and $3\card{\control{G}}-2\chi(G)$ vertices, where $\card{\control{G}}$ is the number of control edges of $G$ and $\chi(G)$ is the Euler characteristic of $G$.

To relate maps of pairs on different spaces, we need the notion of \emph{exact homotopy equivalence} defined in \cite{Collins2001TOPOA}.
A graph map representing the topology of a trellis via exact homotopy equivalence is called \emph{compatible} with the trellis.
\begin{definition}[Compatible graph map]
Let $(G,W)$ be a topological pair where $G$ is a graph and $W$ is a finite subset of $G$.
Then $(G,W)$ is \emph{compatible} with a transverse trellis $T$ if $(G,W)$ and $\cutp{T}$ are exact homotopy equivalent by an embedding $i:(G,W)\exto\cutp{T}$.
A controlled graph map $g$ of $(G,W)$ is \emph{compatible} with the trellis type $[f;T]$ if the embedding $i$ is an exact homotopy equivalence between $g$ and $\cutf{f}$, and $g$ and $\cutf{f}$ have the same orientation at points of $W$.
\end{definition}
Note that the inclusion $i$ induced a bijection between the regions of $\cutt{T}$ (and hence of $T$) with the regions of $(G,W)$,
 and that all compatible graphs are exact homotopy equivalent.
We restrict to transverse trellises since a trellis with tangencies may not have a compatible controlled graph.
While it is always possible to find a topological pair $(G,H)$ which is exact homotopy equivalent to $\cutp{T}$ for which $G$ is a graph,
 it may not be possible to take $H$ to be a finite set of points.

\fig{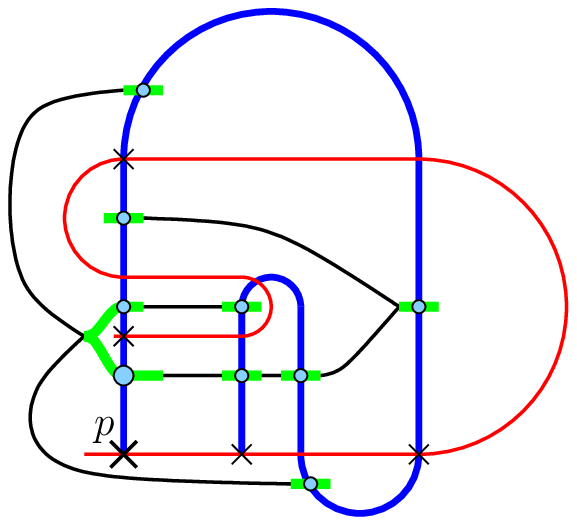}{The controlled graph $(G,W)$ is compatible with the trellis $T$.}
The controlled graph shown in \figref{henonembedded} is compatible with the trellis $T$.
The control edges are shown as thick shaded lines.

There are many controlled graph maps compatible with a trellis type $[f;T]$.
To use graph maps to describe the dynamics forced by $[f;T]$, we need to define a subclass of controlled graph maps which have minimal entropy in the exact homotopy class.
These graph maps are called \emph{efficient}, by analogy with Nielsen-Thurston theory.
Of the efficient graph maps, there is at most one which is \emph{optimal}, giving a canonical \emph{graph representative} for the trellis type.
Tight, efficient and optimal graph maps can be defined in terms of their actions on the turns of $G$.
In this paper, we only need the definition of an optimal graph map; other definitions can be found in \cite{CollinsPPDYNSYS}.
\begin{definition}[Graph representatives]
A controlled graph map $(g;G,W)$ is \emph{optimal} if for every turn $e_1\lhd e_2$, either $\partial g(e_1)\neq \partial g(e_2)$, or at least one of $e_1$ or $e_2$ is a control edge.
A controlled graph map $(g;G,W)$ is a \emph{graph representative} of a transverse trellis type $[f;T]$ if $g$ is an optimal graph map which is compatible with $([f];T)$.
\end{definition}

The following theorem \cite{CollinsPPDYNSYS} shows that every proper irreducible trellis type has a unique graph representative.
\begin{theorem}[Existence and uniqueness of graph representatives]
Let $[f;T]$ be a proper trellis type with no invariant curve reduction.
Then $[f;T]$ has a unique graph representative $(g;G,W)$.
Further, if $[f_0;T_0]$ and $[f_1;T_1]$ are trellis types, the graph representatives $(g_0;G_0,W_0)$ and $(g_1;G_1,W_1)$ are homeomorphic if and only if $[f_0;T_0]=[f_1;T_1]$.
\end{theorem}
In particular, this result shows that the graph representative provides a convenient way of specifying a trellis type.

\pagebreak

%****************************************************************************************************************************************************************************

\section{Regular domains and alpha chains}
\label{sec:regular}

%----------------------------------------------------------------------------------------------------------------------------------------------------------------------------

\subsection{Regular domains}
\label{sec:regulardomains}

Any diffeomorphism is conjugate to a linear hyperbolic diffeomorphism in a neighbourhood of a point of a periodic saddle point.
We therefore expect a trellis map $([f];T)$ to behave in a fairly predictable way in a neighbourhood of a point of $T^P$.
In particular, we know how a rectangle with sides parallel to the local stable and unstable foliations behaves.
Unfortunately, for a given trellis mapping class, the neighbourhood of $T^P$ on which we have hyperbolic behaviour may be arbitrarily small.
To deal with this problem, we introduce the concept of a \emph{regular domain}, which is a rectangular domain which behaves similarly to a sufficiently small rectangular neighbourhood of a periodic saddle point.

\begin{definition}[Regular domain]
\label{defn:regulardomain}
Let $(f;T)$ be a trellis mapping pair.
A rectangular domain $D$ of $T$ is a \emph{regular domain} for a period-$n$ quadrant $Q$ of $T$ if $D$ has sides $T^U[p,q^u]$, $T^S[p,q^s]$, $T^U[q^s,r]$ and $T^S[q^u,r]$, such that 
\begin{enumerate}
\newcounter{uctr}
\newcounter{sctr}
\renewcommand{\theuctr}{(\arabic{uctr}u)}
\renewcommand{\thesctr}{(\arabic{sctr}s)}
\item[(1u)]\refstepcounter{uctr}\label{item:regulardomainadjacentunstable}
  $f^n(T^U[p,q^u])\subset T^U$ with $f^n(T^U[p,q^u])\cap D=T^U[p,q^u]$,
\item[(1s)]\refstepcounter{sctr}\label{item:regulardomainadjacentstable}
  $f^{-1}(T^S[p,q^s])\subset T^S$ with $f^{-n}(T^S[p,q^s])\cap D=T^S[p,q^s]$,
\item[(2u)]\refstepcounter{uctr}\label{item:regulardomainoppositeunstable} 
  $f^{-n}(T^U[q^s,r])\cap D=\emptyset$, and
\item[(2s)]\refstepcounter{sctr}\label{item:regulardomainoppositestable} 
  $f^n(T^S[q^u,r])\cap D=\emptyset$.
\end{enumerate}
\end{definition}
The sides $T^U[p,q^u]$ and $T^S[p,q^s]$ are called \emph{adjacent sides} of $D$, and the sides $T^U[q^s,r]$ and $T^S[q^u,r]$ are \emph{opposite sides}.
The definition is given purely in terms of the topology of the trellis and the mapping of its vertices; no extensions are needed.
Notice that the definition is invariant under time reversal on time reversal if we interchange $T^U$ and $T^S$.
In general, we therefore only need to prove statements on the topology of an adjacent domain for either the unstable (u) or stable (s) case.

\fig{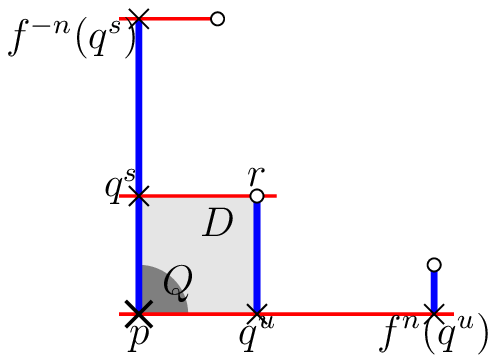}{The domain $D$ is a regular domain for $Q$.}
A regular domain is shown in \figref{regularquadrant}.
We shall always denote the vertices of a regular domain by $p$, $q^u$, $q^s$ and $r$ as shown.
We now give a number of elementary properties of a regular domain.
The first two lemmas are trivial.
\begin{lemma}
\label{lem:regulardomain}
Let $D$ be a regular domain for $(f;T)$ at $Q$. Then
\begin{enumerate}
\item\label{item:regulardomainisotopic} If $f^{-(n+1)}(T^S[p,q^s])\subset T^S$ and $f^{-1}(T^S[q^u,r])$ are subsets of $T^S$, then $f^{-1}(D)$ is a regular domain for $(f;T)$ at $Q$.
\item If $\widetilde{f}$ is isotopic to $f$ relative to $T$, then $D$ is a regular domain for $(\widetilde{f};T)$ at $Q$.
\item If $(\widehat{f};\widehat{T})$ is an extension of $([f];T)$, then $D$ is a regular domain for $(\widehat{f};\widehat{T})$ at $Q$.
\end{enumerate}
\end{lemma}
\begin{proof}
\begin{enumerate}
\item 
Since $f^{-1}(T^U)\subset T^U$, $f^{-1}(T^S[p,q^s])\subset f^{-(n+1)}(T^S[p,q^s])\subset T^S$ and $f^{-1}(T^S[q^u,r])\subset T^S$, the boundary of $f^{-1}(D)$ is composed of segments of $T$.
On the adjacent sides we have $f^n(T^U(f^{-1}(p),f^{-1}(q^u)))\cap f^{-1}(D)=f^{-1}(f^{n}(T^U[p,q^u])\cap D)=f^{-1}(T^U[p,q^u])=T^U[f^{-1}(p),f^{-1}(q^u)]$ and $f^{-n}(T^S[f^{-1}(p),f^{-1}(q^s)]\cap f^{-1}(D))=f^{-1}(f^{-n}(T^U[p,q^s])\cap D)=f^{-1}(T^S[p,q^s])=T^S[f^{-1}(p),f^{-1}(q^s)]$.
On the opposite sides, we have $f^n(f^{-1}(T^S[q^u,r]))\cap f^{-1}(D)=f^{-1}(f^{n}(T^S[q^u,r])\cap D)=\emptyset$ and similarly $f^{-n}(f^{-1}(T^U[q^s,r]))\cap f^{-1}(D)=\emptyset$.
\item
Condition~\ref{item:regulardomainadjacentunstable} shows that $f^n(q^u)\in T^V$, so $\widetilde{f}^n(T^U[p,q^u])=T^U[\widetilde{f}^n(p),\widetilde{f}^n(q^u)]=T^U[f^n(p),f^n(q)]=f^n(T^U[p,q])$, so condition~\ref{item:regulardomainadjacentunstable} is satisfied for $\widetilde{f}$.
If $f^n(r)\in T^V$, then $\widetilde{f}^n(T^S[q^u,r])=f^n(T^S[q^u,r])$, so condition~\ref{item:regulardomainoppositestable} is satisfied.
Otherwise, $f^n(T^S[q^u,r])\subset T^S(f^n(q^u),x)$ for some $x\in T^V$ such that $T^S(f^n(r),x)$ contains no points of $T^U$.
Then since $T^S[f^n(q^u),f^n(r)]\cap D=f^n(T^S[q^u,r])\cap D=\emptyset$, we have $T^S[f^n(q^u),x)\cap D=\emptyset$, and since $\widetilde{f}^n(T^S[q^u,r])\subset T^S[f^n(q^u),x)$, we have $\widetilde{f}^n(T^S[q^u,r])\cap D=\emptyset$, so condition~\ref{item:regulardomainoppositestable} is still satisfied.
A similar argument proves conditions~\ref{item:regulardomainadjacentstable} and \ref{item:regulardomainoppositeunstable}.
\item
By part~\ref{item:regulardomainisotopic}, the domain $D$ is a regular domain for $(\widetilde{f};T)$.
The result follows since the conditions of Definition~\ref{defn:regulardomain} depend only on $T^U$ and $T^S$.
\end{enumerate}
\end{proof}

We now show how to construct regular subdomains of a regular domain.
\begin{lemma}
\label{lem:regularsubdomain}
Let $D$ be a regular domain for $(f;T)$ at a quadrant $Q$.
\begin{enumerate}
\item\label{item:regularsubdomain}
 If $T^S[\tilde{q}^u,\tilde{r}]\subset D$ is a curve with endpoints $\tilde{q}^u\in T^U(p,q^u)$ and $\tilde{r}\in T^U(q^s,r)$, then the rectangular domain $\widetilde{D}$ with vertices at $\{p,\tilde{q}^u,q^s,\tilde{r}\}$ is a regular domain for $Q$.
Similarly, if $T^U[\tilde{q}^s,\tilde{r}]\subset D$ is a curve with endpoints $\tilde{q}^s\in T^S(p,q^s)$ and $\tilde{r}\in T^S(q^u,r)$, then the rectangular domain vertices at $\{p,q^u,\tilde{q}^s,\tilde{r}\}$ is a regular domain for $Q$.
\item\label{item:regulardomainiterate}
 If $f^{-n}(T^S[q^u,r])\subset T^S$, then $f^{-n}(T^S[q^u,r])$ contains a subinterval of the form $T^S[\tilde{q}^u,\tilde{r}]\subset D$ such that $\tilde{q}^u=f^{-n}(q^u)\in T^U(p,q^u)$ and $\tilde{r}\in T^U(q^s,r)$.
Similarly, if $f^{n}(T^U[q^s,r])\subset T^U$, then $f^{n}(T^U[q^s,r])$ contains a subinterval of the form $T^U[\tilde{q}^s,\tilde{r}]\subset D$ such that $\tilde{q}^s=f^{n}(q^s)\in T^S(p,q^s)$ and $\tilde{r}\in T^S(q^u,r)$.
\end{enumerate}
\end{lemma}
\begin{proof}
\begin{enumerate}
\item 
The only nontrivial step is to show $f^n(T^S[\tilde{q}^u,\tilde{r}^u])\cap\widetilde{D}=\emptyset$.
First, note that $f^n(\tilde{q}^u)\in T^U(\tilde{q}^u,f^n(q^u))$, so $f^n(\tilde{q}^u)\not\in\widetilde{D}$.
Further, $f^n(T^S[\tilde{q}^u,\tilde{r}^u])\cap T^U[p,\tilde{q}^u]=f^n(T^S[\tilde{q}^u,\tilde{r}^u]\cap f^{-n}(T^U[p,\tilde{q}^u]))=f^n(\emptyset)=\emptyset$
 and $f^n(T^S[\tilde{q}^u,\tilde{r}^u])\cap T^U[p,\tilde{q}^u]=f^n(T^S[\tilde{q}^u,\tilde{r}^u]\cap f^{-n}(T^U[p,\tilde{q}^u]))\subset f^n(D\cap f^{-n}(T^U[p,\tilde{q}^u))=\emptyset$, so $f^n(T^S[\tilde{q}^u,\tilde{r}^u])$ does not intersect the unstable boundary of $\widetilde{D}$.
Additionally, $f^n(T^S[\tilde{q}^u,\tilde{r}])\cap T^S[p,q^s]=f^n(T^s[\tilde{q}^u,\tilde{r}]\cap f^{-n}(T^S[p,q^s]))\subset f^n(T^S[\tilde{q}^u,\tilde{r}]\cap T^S[p,q^s])=\emptyset$.
It remains to show that $f^n(T^S[\tilde{q}^u,\tilde{r}])$ does not intersect itself.
We have already seen that $f^n(T^S[\tilde{q}^u,\tilde{r}])$ does not contain $\tilde{q}^u$ or $\tilde{r}$.
Further, $f^n(\{\tilde{q}^u\})\cap T^S[\tilde{q}^u,\tilde{r}]\subset f^n(T^U(p,q^u))\cap T^S[\tilde{q}^u,\tilde{r}]=\{\tilde{q}\}$, but clearly $f^n(\tilde{q}^u)\neq\tilde{q}^u)$, so $f^n(\tilde{q}^u)\not\in T^S[\tilde{q}^u,\tilde{r}]$.
Since $T^S[\tilde{q}^u,\tilde{r}]$ is an interval, this is enough to show that $f^n(T^S[\tilde{q}^u,\tilde{r}])\cap T^S[\tilde{q}^u,\tilde{r}]$.
Hence $f^n(T^S[\tilde{q}^u,\tilde{r}])\cap \partial\widetilde{D}=\emptyset$, so $f^n(T^S[\tilde{q}^u,\tilde{r}])\cap \widetilde{D}=\emptyset$.
A similar analysis proves the statement for a curve $T^U[\tilde{q}^s,\tilde{r}]$.
\item
Since $f^n(Q)=Q$, $f^n$ is orientation-preserving, so the orientation of the intersection at $f^{-n}(q^u)$ is the same as that at $q^u$.
Further, $f^{-n}(T^S[q^u,r])\cap T^U[p,q^u]=f^{-n}(T^S[q^u,r]\cap f^n(T^U[p,q^u]))=f^{-n}(T^S[q^u,r]\cap T^U[p,q^u])=\{f^{-n}(q^u)\}$, so $f^{-n}(T^S(q^s,r])$ does not cross $T^U[p,q^u]$.
Finally, $\{f^{-n}(r)\}\cap D\subset f^{-n}(T^U[q^s,r])\cap D=\emptyset$, so $f^{-n}(r)\not\in D$.
Therefore, $f^{-n}(T^S[q^u,r])$ intersects $T^U[q^s,r]$, and we let $\tilde{r}$ be this first intersection.
\end{enumerate}
\end{proof}

\fig{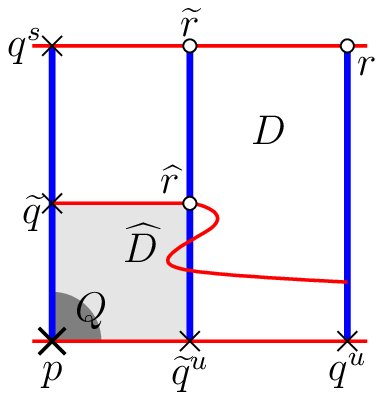}{$\widehat{D}$ is a regular domain if $D$ is a regular domain.}
In Figure~\ref{fig:regularsubdomain} we show a regular domain $D$ of $(f;T)$ at a quadrant $Q$. 
An application of Lemma~\ref{lem:regularsubdomain} shows that the rectangular domain with vertices at $p$, $\widetilde{q}^u$, $q^s$ and $\widetilde{r}$ is a regular domain for $Q$.
A further applications of Lemma~\ref{lem:regularsubdomain}(\ref{item:regularsubdomain}) shows that the domain $\widetilde{\widetilde{D}}\subset D$ is also a regular domain.
By continuing such a construction, it is possible to find a regular domain for $Q$ contained in an arbitrarily small neighbourhood of $p$.

\fig{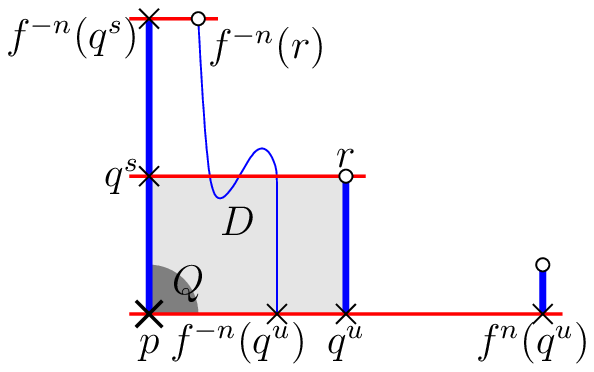}{The backward iterate of $T^S[q^u,r]$ intersects $T^U[q^s,r]$ but does not cross $T^U[p,q^u]$.}
In \figref{regularsideiterate} we show a regular domain $D$ together with $f^{-n}(T^S[q^u,r])$.
Notice that $f^{-n}(T^S[q^u,r])$ may cross $T^U[q^s,r]$ several times, but must cross at least once since $f^{-n}(r)\not\in D$.

\begin{remark}
Lemma~\ref{lem:regularsubdomain}(\ref{item:regulardomainiterate}) is enough to show that if $D$ is a regular domain for a trellis mapping class $([f];T)$, and $([\widetilde{f}];\widetilde{T})$ is an extension of $([f];T)$ with $\widehat{f}^{-n}(T^S[q^u,r])\subset \widehat{T}^S$, then the biasymptotic orbit through $\tilde{r}$ is forced by $([f];T)$.
This result has potential applications in the theory of forcing relations for homoclinic and periodic orbits \cite{deCarvalhoHall2002EXPMAT}.
\end{remark}

%----------------------------------------------------------------------------------------------------------------------------------------------------------------------------

\subsection{Alpha-chains and transitivity}

The most important case of a regular domain for a quadrant $Q$ is when $D$ is a region, in which case we have a \emph{regular region} for $Q$.
It is also useful to consider the less restrictive case for which the interior of $D$ does not intersect $T^U$, or in other words, if $D\cap T^U\subset\partial D$.
It is an important but trivial observation that if $D\cap T^U=\partial^D D$ and if $\widehat{T}$ is an extension of $T$ with $\widehat{T}^U=T^U$, then $D\cap\widehat{T}^U=\partial^U D$.

If $D$ is a regular domain such that $D\cap T^U=\partial^U D$
In this case, there is a unique homotopy class of exact curves $\alpha_D:(I,J)\exto\cutp{T}$ such that $\alpha_D(I)\subset D$, $\alpha_D(0)\in T^U[p,q^s]$, $\alpha_D(1)\in T^U[q^u,r]$, and $\alpha_D$ has minimal intersections with $T^S$.
The first intersection of $\alpha_D$ with $T^S$ must be in a segment of $T^S$ from $T^U[p,q^u]$ to $T^U[q^s,r]$, and this segment forms the opposite stable side $T^S[\widetilde{q}^u,\widetilde{r}]$ of a regular domain $\widetilde{D}\subset D$.
For any regular domain $D$ for $(f;T)$ at a period-$n$ quadrant $Q$, the $n$th minimal iterate $f_{\min}^n[\alpha_D]$ of $[\alpha_D]$ tightens onto $[\alpha_D]$ showing that the minimal iterates of $[\alpha_D]$ expand.
If $f_{\min}^n[\alpha_D]$ tightens to $[\alpha_{D\p}]$ for some other segment $D\p$, we can deduce information on the crossing of $f^k(D)$ with $D\p$.
\fig{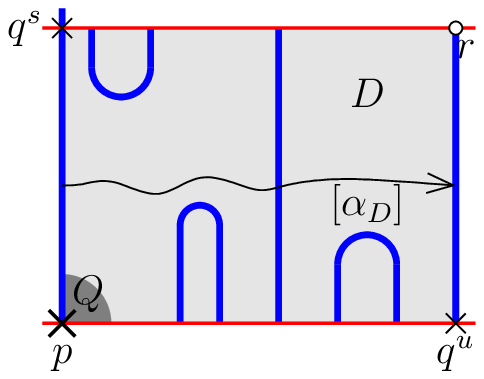}{The homotopy class $[\alpha_D]$.}

\begin{definition}[Minimal regular domain crossing curve]
Let $D$ be a regular domain for a periodic-$n$ quadrant $Q$ under $(f;T)$. 
Then the homotopy class $[\alpha_D]$ is the homotopy class of an exact curve $\alpha_D:(I,J)\fto \cutt{T}$ from $T^S[p,q^s]$ to $T^S[q^u,r]$ in $D$ which has minimal intersections with $T^S$.
\end{definition}

If $Q$ is a quadrant with a regular domain $D$ such that $D\cap T^U=\partial^U D$, we let $D(Q)$ be the smallest regular domain for $Q$, and write $[\alpha_Q]$ for $[\alpha_{D(Q)}]$.
Since for any regular domain $D$, the $n$th iterate of $\alpha_D$ tightens onto itself, there is some recurrent behaviour for the graph representative.

The following lemma shows that any curve in such a regular domain with initial point in $T^S[p,q_s]$ and and endpoint in some other interval of $T^S$ iterates to a curve which crosses $T^S[q^u,r]$.
We then show that if $D$ and $\widetilde{D}$ are regular domains for a quadrant, then some iterate of $\alpha_D$ tightens onto $\alpha_{\widetilde{D}}$.
The importance of this result is that if a curve $\widehat{\alpha}$ tightens onto $\alpha$, and both curves are tight curves embedded in a graph, then $\alpha(I)\subset\widehat{\alpha}(I)$.
\begin{lemma}\label{lem:alphacurve}
\ 
\begin{enumerate}
\item\label{item:subdomainiterate}
Let $D$ be a regular domain for a quadrant $Q$ such that $D\cap T^U=\partial^U D$ and $D$ contains no smaller regular domains, and $D\cap T^P=\{p\}$.
Then for any exact curve $\alpha:(I,\{0,1\})\fto\cutp{T}$ with $\alpha(I)\subset D$, $\alpha(0)\in T^S[p,q^s]$ and $\alpha(1)\not\in T^S[p,q^s]$, there exists $k$ such that $f_{\min}^{nk}[\alpha]$ crosses $T^S[q^u,r]$.
\item\label{item:subdomaintighten}
Let $D$ be a regular domain for a quadrant $Q$ such that $D\cap T^U=\partial^U D$, and let $\widetilde{D}\subset D$ be a regular domain.
Then there exists $k$ such that $f_{\min}^{kn}[\alpha_{D}]$ tightens onto $[\alpha_{\widetilde{D}}]$.
\item\label{item:subdomainiteratetighten}
Let $D$ be a regular domain for a quadrant $Q$, and $\widetilde{D}\subset D$ a regular domain for $Q$ such that $\widetilde{D}\cap T^U=\partial\widetilde{D}^U$.
Let $\alpha$ be a curve in $D$ such that $\alpha(0)\in T^S[p,q^s]$ and $\alpha(1)\in T^S[q^u,r]$.
Then there exists $k$ such that $f_{\min}^{kn}[\alpha]$ tightens onto $[\alpha_{\widetilde{D}}]$.
\end{enumerate}
\end{lemma}
\begin{proof}
\begin{enumerate}
\item
Let $S$ be the segment containing the final endpoint of $\alpha$.
Suppose $f^{in}(S)\subset D$ for all $i$.
Then $f^{in}(S)\cap T^S[p,q^s]=\emptyset$ for all $i$, and since there are only finitely many stable segments in $D$, one must map into itself under $f^{jn}$ for some $j$, giving a periodic point in $D$ distinct from $p$, a contradiction.
Therefore $f^{kn}(S)\not\subset D$ for some least $k$
For this $k$, we must have $f_{\min}^{kn}[\alpha]\cap T^S[q^u,r]\neq\emptyset$, since crossing $T^S[p,q^s]$ would imply a point of $f^{kn}\circ\alpha$ in $T^S(q^s,f^{-n}(q^s))$, a contradiction.
\item
Let $k$ be the least integer such that $f^{kn}(\tilde{q}^u)\not\in T^U[p,q^u]$.
Then $f^{kn}(\tilde{q}^u)\not\in D$, so the endpoint of $f_{\min}^{kn}[\alpha_{\widetilde{D}}]$ does not lie in $D$.
Further, the first intersection of $f_{\min}^{kn}[\alpha_{\widetilde{D}}]$ with $\partial^S D$ must be with $T^S[q^u,r]$, so $f_{\min}^{kn}[\alpha_{\widetilde{D}}]$ tightens onto $[\alpha_D]$, which is the curve from $T^S[p,q^s]$ to $T^S[q^u,r]$ in $D$ with minimal intersections with $T^S$.
\item
For any $i$, $f_{\min}^{in}[\alpha]$ must cross $T^S[q^u,r]$, and do so before it crosses $T^S[p,q^s]$.
Choose $k$ such that $f_{\min}^{kn}(\alpha(0))\in T^S[p,\tilde{q}^s]$.
Then $f_{\min}^{kn}[\alpha]$ tightens onto $[\alpha_{\widetilde{D}}]$.
\end{enumerate}
\end{proof}

\fig{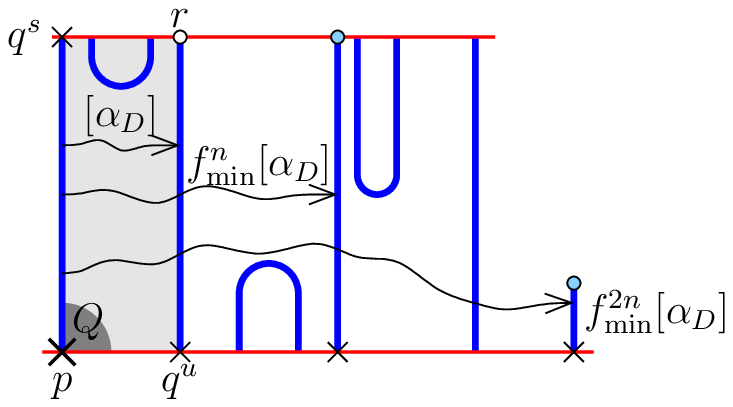}{A curve $\alpha_D$ in a regular domain $D$ and some of its minimal iterates.}
An example of the minimal iterates of a curve $\alpha=\alpha_D$ for a regular domain $D$ is shown in \figref{alphacover}.

If a homotopy class $[\alpha_D]$ has some iterate which tightens onto $\alpha_D$ for a different quadrant, the graph representative maps edges from one quadrant to another.
This provides the definition of an \emph{alpha chain}.
\begin{definition}[Alpha chain]
A sequence of homotopy classes of exact curves $[\alpha_i]$ for $0\leq i\leq k$ with endpoints in $T^S$ is an \emph{alpha chain} if $f_{\min}[\alpha_i]$ tightens onto $[\alpha_{i+1}]$ for $0\leq i<k$.
We say there is an \emph{alpha chain} from quadrant $Q_U$ to quadrant $Q_S$ under the isotopy class $([f];T)$ if $Q_U$ and $Q_S$ are contained in regular domains $D_U$ and $D_S$ which do not intersect $T^U$ in their interiors, and there is an alpha chain from $[\alpha_{D_U}]$ to $[\alpha_{D_S}]$.
\end{definition}
\fig{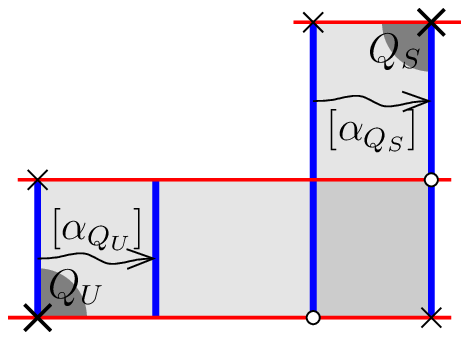}{A trellis with an alpha chain from $Q_U$ to $Q_S$.}
An example of an alpha chain from a quadrant $Q_U$ to $Q_S$ is shown in \figref{alphachain}.
The following lemma gives a sufficient condition for the existence of an alpha-chain from $Q_U$ to $Q_S$.
\begin{lemma}
Let $D_U$ and $D_S$ be regular domains for $Q_U$ and $Q_S$ respectively.
\begin{enumerate}
\item
If $D_U\cap D_S$ is a rectangular domain with unstable edges contained in $\partial^U D_U$ and stable edges in $\partial^S D_S$, then there is an alpha chain from $Q_U$ to $Q_S$.
\item
If there is an alpha chain from $Q_U$ to $Q_S$ for $([f];T)$ and $([\widehat{f}];\widehat{T})$ is a stable extension of $T$, then there is an alpha chain from $Q_U$ to $Q_S$ for $([\widehat{f}];\widehat{T})$.
\end{enumerate}
\end{lemma}
\begin{proof}
\begin{enumerate}
\item
Let $\beta$ be a curve in $D_U\cap D_S$ homotopic to the initial interval of $T^U(p_u)\cap D_S$.
Then $f_{\min}^{k_un_u}[\alpha_{Q_U}]$ tightens onto $\beta$ for some $k_u$, and since $\beta$ is a curve joining the stable sides of $D_S$, we also have $f_{\min}^{k_sn_s}[\beta]$ tightens onto $\alpha_{Q_S}$ for some $k_s$.
Hence $f_{\min}^{k_un_u+k_sn_s}[\alpha_{Q_U}]$ tightens onto $[\alpha_{Q_S}]$ as required.
\item
By Lemma~\ref{lem:alphacurve}(\ref{item:subdomaintighten}), $\widehat{f}_{\min}^{k_un_u}[\widehat{\alpha}_{Q_U}]$ tightens onto $[\alpha_{Q_U}]$ for some $k_u$, and by Lemma~\ref{lem:alphacurve}(\ref{item:subdomainiteratetighten}), $\widehat{f}_{\min}^{k_sn_s}[\alpha_{Q_S}]$ tightens onto $[\alpha_{Q_S}]$ for some $k_s$.
Further, $f_{\min}^{n}[\alpha_{Q_U}]$ tightens onto $[\alpha_{Q_S}]$ for some $n$ since there is an alpha-chain from $Q_S$ to $Q_S$ in $T$.
Therefore, there is an alpha chain from $Q_U$ to $Q_S$ in $\widetilde{T}$.
\end{enumerate}
\end{proof}

\begin{definition}[Transitive trellis]
We say that a trellis type $[f;T]$ is \emph{transitive} if every quadrant $Q$ is contained in a regular domain regular $D_Q$ such that $D_Q\cap T^U=\partial D_Q$, and for every pair of quadrants $Q_U$ and $Q_S$, there is an alpha chain from $Q_U$ to $Q_S$.
\end{definition}

Transitivity of a trellis type is an especially useful property since it means that the trellis type has a transitive graph representative.
We now show that the hyperbolicity near $T^P$ is enough to create intersections from which we can deduce regularity and transitivity.
Note that for this result we are concerned with $f$-extensions of $(f;T)$ rather than an isotopy class.
\begin{lemma}
\label{lem:regularextension}
Let $T$ be a trellis for a diffeomorphism $f$, and let $Q$ a quadrant of $T$.
Let $q\in T^U(Q)$ be the endpoint of an unstable segment $S$ on the $Q$-side of $T^S(Q)$, and let $q^s\in T^S(Q)$ be the endpoint of a stable segment $U$ on the $Q$-side of $T^S(Q)$.
Then there exists $k$ such that $f^{-kn}(S)$ intersects $U$ at a point $r$ such that $\{p,f^{-kn}(q),q^s,r\}$ are the vertices of a regular domain for $Q$.
\end{lemma}

\begin{proof}
By the Lambda lemma, as $i\tendsto\infty$, $f^{-in}(S)$ limits on $W^S(Q)$ in the $C^1$ topology.
Take a neighbourhood $K$ of $q^s$ such that $f^{-n}(K)\cap T^S[p,q^s]=\emptyset$, and choose $k$ such that $f^{-kn}(S)$ intersects $U$ in at a point $r$ in $K$ such that $T^S[q^s,r]\subset K$, and the domain $D$ with vertices at $\{p,f^{-kn}(q),q^s,r\}$ is a rectangle which does not intersect $f^{-n}(K)$, and such that $f^n(T^U(p,f^{-kn}(q^u)))\cap D=T^U[p,f^{-kn}(q^u)]$.
Let $\widehat{T}=(T^U,f^{-nk}(T^S))$.
Then $f^{-n}(T^U[q^s,r])\cap D\subset f^{-n}(K)\cap D=\emptyset$, so $D$ is a regular domain of $(f;\widehat{T})$, as shown in \figref{regularisequadrant}.
\fig{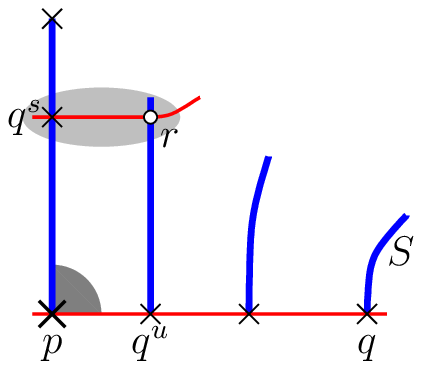}{Backward iterates of $S$ eventually form the opposite stable side of a regular domain.}
\end{proof}

\begin{corollary}
\label{cor:regularextension}
If $(f;T)$ is irreducible, then there is an $f$-extension $\widehat{T}$ of $T$ such that $(f;\widehat{T})$ is transitive.
\end{corollary}

%****************************************************************************************************************************************************************************

\section{Entropy-Minimising Diffeomorphisms}
\label{sec:entropy}

In this section, we show that the entropy bound obtained by the Nielsen entropy is sharp, at least for well-formed trellises.
That is, the topological entropy of the graph representative $g$ for a trellis type $[f;T]$ (which is the same as the Nielsen entropy $\hniel[f;T]$)
 is the infemum of the topological entropies of diffeomorphisms in the class.
The main results of the paper are summarised in the following theorem.

\begin{theorem*}
\label{thm:mainentropy}
Let $([f];T)$ be a well-formed trellis type.
Then for any $\epsilon>0$ there exists a diffeomorphism $\widehat{f}\in[f;T]$ such that $\htop(\widehat{f})<\hniel[f;T]+\epsilon$.
If there exists a diffeomorphism $\widehat{f}$ isotopic to $f$ relative to $T$ such that any $\widehat{f}$-extension of $T$ is minimal, then there exists a uniformly-hyperbolic diffeomorphism $\widetilde{f}\in[f;T]$ such that $\htop(\widetilde{f})=\hniel[f;T]$.
Further, if $[f;T]$ is irreducible and $\widetilde{f}\in[f;T]$ such that $\htop(\widetilde{f})=\hniel[f;T]$, then any $\widetilde{f}$-extension of $T$ is minimal.
\end{theorem*}

We now give some examples which illustrate the hypotheses of the theorem.
The following example shows that the hypothesis that the trellis be well-formed is necessary.
\begin{example}
\label{ex:illformed}
The trellis mapping classes in \figref{illformed} are not well-formed.
\fig{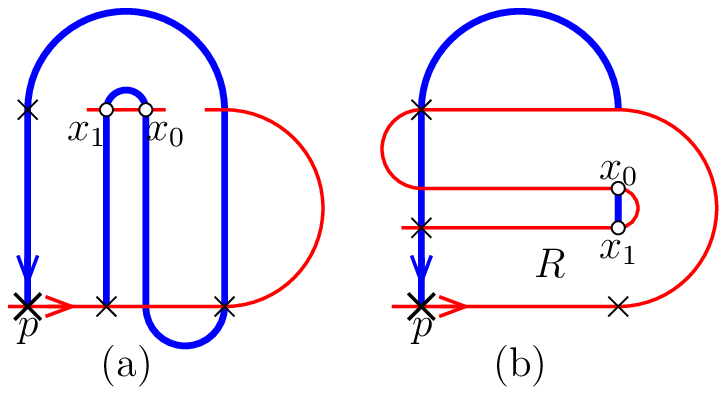}{Two ill-formed trellis mapping classes. The trellis in (a) has Nielsen entropy $\log2$, whereas the trellis in (b), which is the time reversal of that in (a), has Nielsen entropy zero.}
\par
The Nielsen entropy of the trellis mapping class in (a) is equal to $\log2$, so any diffeomorphism in the trellis mapping class must have topological entropy at least $\log2$.
Since the Smale horseshoe map has this trellis type, the topological entropy of the trellis type is exactly $\log2$.
The trellis mapping class of \figrefpart{illformed}{b} is conjugate to the time-reversal of the trellis mapping class in (a).
Since the topological entropy of a diffeomorphism is the same as that of its inverse, any diffeomorphism in this trellis mapping class must have topological entropy at least $\log2$.
However, all the edges of the graph representative are control edges, so the Nielsen entropy is zero.
\end{example}
The above example illustrates that a trellis which is not well-formed may have Nielsen entropy strictly less than the topological entropy, and may even have different Nielsen entropy from its time-reversal.

Even if a trellis mapping class is well-formed, it is not necessarily true that the Nielsen entropy is realised.
A trivial example is of a well-formed trellis type $[f;T]$ for which $\hniel[f;T]=\htop[f;T]$ by for which every diffeomorphism in $[f;T]$ has topological entropy greater than $\hniel[f;T]$ is the planar trellis type with a single transverse homoclinic intersection.
The Nielsen entropy of this trellis type is equal to zero, but every diffeomorphism with a transverse homoclinic point has strictly positive topological entropy.
However, it is simple to construct trellis maps with topological entropy arbitrarily close to zero. 
We now give a nontrivial example.
\begin{example}
\label{ex:nominimalentropy}
\fig{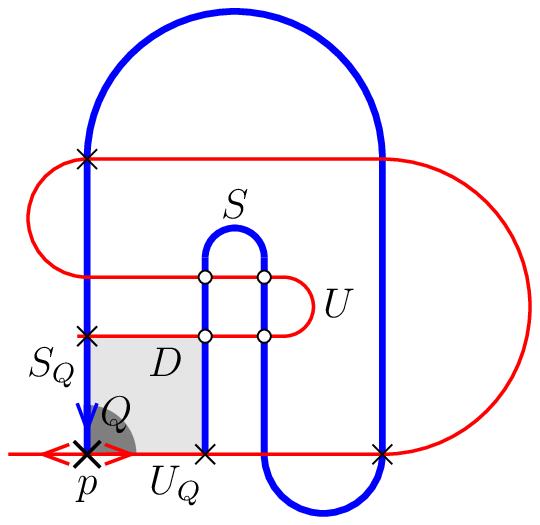}{A trellis for which the Nielsen entropy is not realisable.}
A trellis type $[f;T]$ for which the Nielsen entropy is not realisable is shown in \figref{nominimalentropy}.
Consider the segment $S$.
Taking backward minimal iterates of $S$ eventually yields a segment lying in the region $D$ with endpoints in the segment $U(Q)$.
Then, by the Lambda lemma, under any diffeomorphism $f$ in the trellis mapping class, $f^{-n}(S)$ tends to the closed branch of $T^S(p)$ as $n\tendsto-\infty$, so contains an intersection with $T^U$ for some $n$, even though any \emph{minimal} backward iterate of $S$ has no intersections with $T^U$.
Similarly, $f^n(U)$ must intersect $T^S$ for some $n$ even though any minimal iterate does not.
These observations can be used to show that the Nielsen entropy of the extension $[f;\widehat{T}]$ is greater than that of $[f;T]$, so $\htop(f)>\hniel[f;T]$
\end{example}

In general, the realisability of the entropy bound is closely related to the existence of a diffeomorphism for which every extension is a minimal extension.
In the case where this infemum is realised, we show how to construct a minimal-entropy uniformly-hyperbolic diffeomorphism in the trellis mapping class.
Otherwise, we show, for any $\epsilon>0$, how to construct a diffeomorphism whose entropy is within $\epsilon$ of the Nielsen entropy.

%----------------------------------------------------------------------------------------------------------------------------------------------------------------------------

\subsection{Construction of minimal supertrellises}
\label{sec:minimalsupertrellis}

Most of the procedures we use to construct diffeomorphisms in a given trellis mapping class rely on extending the original trellis
 and introducing new branches in a controlled way.
The most important type of supertrellis is a minimal supertrellis, since we expect the Nielsen entropy of a minimal supertrellis to be the same as that of the original trellis mapping class.
This is indeed the case, as the following theorem from \cite{CollinsPPDYNSYS} shows.
\begin{theorem}[Nielsen entropy of minimal supertrellises]
\label{thm:minimalsupertrellisentropy}
Let $([f];T)$ be a well-formed trellis mapping class.
If $([\widetilde{f}];\widetilde{T})$ is a minimal supertrellis of $([f];T)$, then $\hniel[\widetilde{f};\widetilde{T}]=\hniel[f;T]$.
\end{theorem}
We then show that we can introduce new stable and unstable branches at essential Nielsen classes of $f$.
The only difficulty here is on finding the correct initial segment of a branch; once this has been achieved, we can take minimal iterates in the usual way.
Note that a stable supertrellis which is not minimal may change the topology of the graph.

\begin{lemma}
\label{lem:invariantarcs}
Let $\{\alpha_i:i=0\ldots n-1\}$ be a set of curves with endpoints in $f(T^S)$ such that each $\alpha_i$ has minimal intersections with $T^S$ (taken homotopies through curves with endpoints in $T^S$) as curves relative to endpoints, and $f_{\min}[\alpha_i]$ tightens onto $[\alpha_{i+1}]$ (using arithmetic modulo $n$).
Suppose further that any $\alpha_i$ with an endpoint in a segment of $T^S$ containing a point of $T^P$ contains a branch of $T^U$.
Then $f$ is isotopic relative to $T$ to a diffeomorphism $\widetilde{f}$ with $\widetilde{f}\circ\alpha_i\supset\alpha_{i+1}$.
The trellis mapping class $([\widetilde{f}];\widetilde{T})$ with $\widetilde{T}=(T^U\cup\bigcup_{i}\alpha_i(I),T^S)$ is a minimal supertrellis of $([f];T)$.
\end{lemma}
\begin{proof}
Let $\beta_{i+1}\in f_{\min}[\alpha_i]$ and be disjoint from $T^U$.
Since $f_{\min}[\alpha_i]$ tightens onto $\alpha_{i+1}$, without loss of generality we can assume that $\beta_{i+1}\supset \alpha_{i+1}$ for all $i$.
The curves $f\circ\alpha_i$ and $\beta_{i+1}$ are isotopic by an isotopy through curves with endpoints in $f(T^S)$.
By the isotopy extension theorem this isotopy can be extended to an ambient isotopy $h_t$ which is fixed on $T^U$ and maps $f(T^S)$ into $f(T^S)$.
Further, since the curves $\alpha_i$ either contain a branch of $T^U$ or lie in segments of $T^S$ which do not contain points of $T^P$, we can ensure that $h_t$ is fixed on the segments of $T^S$ containing points of $T^P$.
Let $\widetilde{f}=h_1\circ f$, so $T^U$ and $T^S$ is are still part of the unstable and stable manifolds of $T^P$ for $\widetilde{f}$.
Then $\widetilde{f}(T^U)=h_1(f(T^U))\supset h_1(T^U)=T^U$, $\widetilde{f}(T^S)=h_1(f(T^S))=f(T^S)\subset T^S$ and $\widetilde{f}(\alpha_i)=h_1(f(\alpha_i))=\beta_{i+1}\supset \alpha_{i+1}$.
Therefore $([\widetilde{f}];\widetilde{T})$, where $\widetilde{T}=(T^U\cup\bigcup\beta_i,T^S)$, is a well-formed trellis mapping class, and is a minimal supertrellis of $([f];T)$ since the curves $\beta_{i+1}$ are minimal iterates of the curves $\alpha_{i}$.
\end{proof}
By reversing time, we obtain a similar result for curves with endpoints in $T^U$ which are disjoint from $T^S$.
Such curves lift to cross cuts in $\cutm{T}$.

\fig{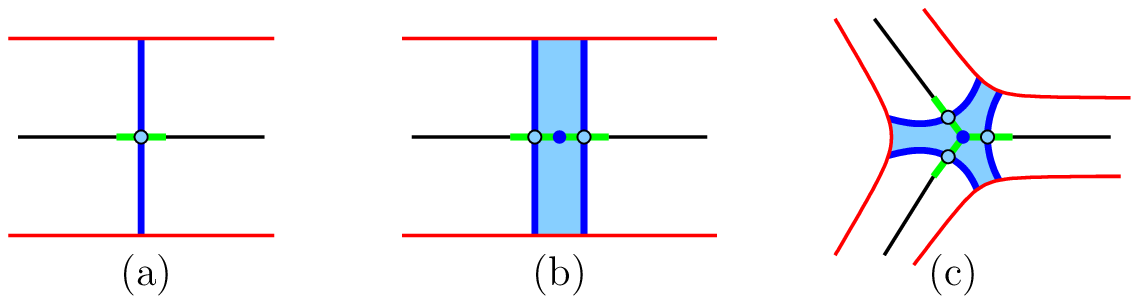}{Trellises and graphs formed by introducing new stable curves. The trellis in (a) is formed by introducing a single stable curve and does not have any attractors. The trellises in (b) and (c) are formed by introducing one stable curve crossing each incident edge at a vertex, and do yield new attractors.}

We can apply the above result in the following concrete cases which we shall need later.
\begin{theorem}
\label{thm:minimalsupertrellis}
Let $([f];T)$ be a well-formed irreducible trellis mapping class.
Then:
\begin{enumerate}
\renewcommand{\theenumi}{(\alph{enumi})}
\item\label{item:periodicblowup}
If $p$ is an essential periodic point of $([f];T)$ which does not shadow $T^P$, there is a minimal supertrellis of $([f];T)$ for which $p$ is contained in a region which can be chosen to be attracting, repelling, or both.
\item\label{item:periodictrivialend} If $T^{\us}[p,b]$ is a trivial branch of $T$, then there is a then there is a minimal supertrellis $([\widetilde{f}];\widetilde{T})$ of $([f];T)$ for which $T^{\us}[p,b]$ is contained in a non-chaotic region of $([\widetilde{f}];\widetilde{T})$.
\item\label{item:trivialbranch} If $T^{\us}[p,b]$ is a trivial branch of $T$ which lies in a chaotic region of $T$, then there is a minimal supertrellis $([\widetilde{f}];\widetilde{T})$ of $([f];T)$ for which $\widetilde{T}^{\us}(p,b)$ is contained in a nontrivial branch of $T^P$.
\end{enumerate}
\end{theorem}
\begin{proof}
In all cases, we suppose that $p$ is contained in a region $R$ and has period $n$.
We only consider the construction of a minimal supertrellis with new unstable curves; the construction of new unstable curves follows by reversing time.
We say a stable segment $S$ on the boundary of $R$ is an \emph{exit segment} for $p$ if the homotopy class of curves $\gamma$ from $p$ to $S$ is such that $f_{\min}^{nk}[\gamma]$ has initial arc $[\gamma]$, but is not equal to $[\gamma]$.
Note that if $[\gamma]$ is a homotopy class of curves from $p$ to an exit segment $S$, then the first intersection of $f_{\min}^{i}[\gamma]$ is an exit segment for $f^i(p)$ for any $i$.
If there are no exit segments for $p$, then $p$ is contained in a region $R$ which is attracting, and hence is non-chaotic.
The results follow by constructing curves the $\alpha_i$ from curves $\gamma$ joining $p$ to exit segments, and applying Lemma~\ref{lem:invariantarcs}.
\begin{enumerate}
\renewcommand{\theenumi}{(\alph{enumi})}
\item
Since $p$ does not shadow $T^P$ and is essential, it is easy to see that there must be at least two exit segments for $p$.
Let $[\gamma_0]$ be a curve joining $p$ to an exit segment, and let $[\gamma_i]$ be the homotopy class of initial arc of $f_{\min}^{i}[\gamma_0]$.
We take curves $\alpha_{i}$ to be the boundary curves of a small neighbourhood of the curves $\gamma_{i}$ for $0\leq i<nk$.
\item
Suppose $p$ lies at the end of a trivial stable branch of $p$.
If there is a single exit segment for $p$, then $p$ already lies in a repelling region.
Otherwise, we can find a curve $\alpha$ joining two exit segments such that the component of $R\setminus \alpha_0(I)$ containing $p$ is as small as possible.
The curves $\alpha_i$ are constructed by iterating $\alpha$.
\item Construct $\alpha_i$ by iterating a curve $\gamma$ from $p$ be an exit segment.
\end{enumerate}
\end{proof}

\fig{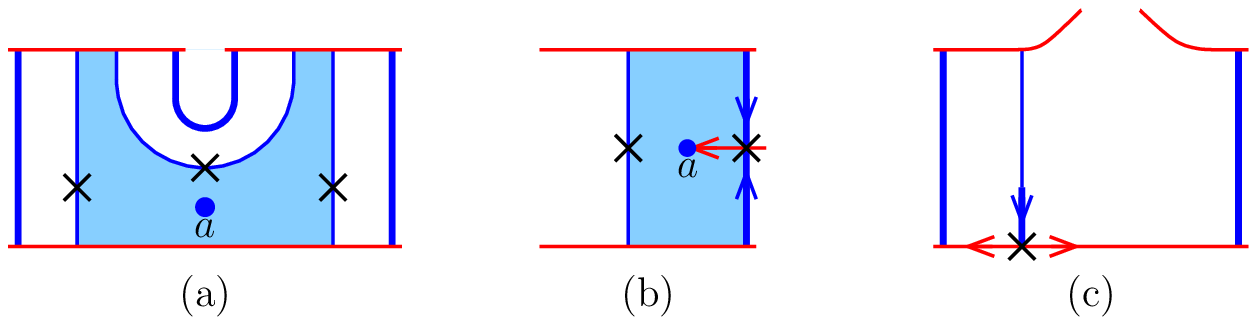}{Stable supertrellises. In (a) we introduce an attractor, in (b) we bound an attractor, and in (c) we extend a trivial branch.}
We use these results to construct new branches at essential periodic orbits and trivial branches of a trellis.
These cases are shown in \figref{stablesupertrellis}.
These branches are may used to create new attractors and repellors, or subdivide a region into an attractor/repellor and a collection of rectangles, as in \figrefpart{stablesupertrellis}{a}, to bound the end of a trivial branch in an attractor or repellor, as in (b), or to extend a trivial branch to a nontrivial branch, as in (c).

%----------------------------------------------------------------------------------------------------------------------------------------------------------------------------

\subsection{Existence of entropy minimisers}

To prove the existence of a diffeomorphism in a trellis mapping class whose topological entropy is the Nielsen entropy of the class, we reduce to the case for which every chaotic region is a rectangle.
We then construct such a diffeomorphism for a particularly simple class of trellises.
\begin{theorem}[Existence of entropy minimisers]
\label{thm:minimalmodel}
Let $([f];T)$ be a trellis mapping class.
Suppose there is a diffeomorphism $\widehat{f}\in([f];T)$ such that every extension of $T$ by $\widehat{f}$ is minimal.
Then there is a uniformly-hyperbolic diffeomorphism $\widetilde{f}\in([f];T)$ such that every extension of $T$ by $\widetilde{f}$ is minimal, and $\htop(\widetilde{f})=\hniel[f;T]$.
\end{theorem}

\begin{proof}
By Theorem~\ref{thm:minimalsupertrellisentropy}, any trellis mapping class of the form $[\widehat{f};(T^U,\widehat{f}^{-n}(T^S)]$ has the same Nielsen entropy at $([f];T)$.
By Lemma~\ref{lem:regularextension} we can therefore take a $\widehat{f}$-extension $\widehat{T}_1$ of $T$ such that every quadrant of $\widehat{T}_1$ lies in a regular domain, 
 and by irreducibility, we can ensure that every region of  $\widehat{T}_1$ is a topological disc or annulus.
If $r_0$ is a point at the end of a trivial branch of $T^S$ and lies in a chaotic region, by introducing new unstable curves as in Theorem~\ref{thm:minimalsupertrellis}\ref{item:periodictrivialend}, we can take a minimal supertrellis $([\widetilde{f}_2];\widetilde{T}_2)$ such that the orbit of $r_0$ lies in a repelling region.
Similarly, if $a_0$ is a point of $\widetilde{T}_2^S$ at the end of a trivial branch of $\widehat{T}_1^U$, by introducing new stable curves, we can take a minimal supertrellis $([\widetilde{f}_3];\widetilde{T}_3)$ such that the orbit of $a_0$ lies in a attracting region.

Now suppose there is a chaotic region $R$ of $([\widetilde{f}_3];\widetilde{T}_3)$ which is not a rectangle.
Then the graph representative $(\widetilde{g}_3];\widetilde{G}_3,\widetilde{W}_3)$ of $([\widetilde{f}_3];\widetilde{T}_3)$ has a peripheral loop or a valence-$n$ vertex in $R$ which corresponds to a boundary component or essential periodic orbit which does not shadow $T^S$.
Introducing new stable curves for all such $R$ as in Theorem~\ref{thm:minimalsupertrellis}\ref{item:periodicblowup} gives a minimal supertrellis $([\widetilde{f}_4];\widetilde{T}_4)$ for which $R$ is a domain containing an attractor and some chaotic rectangles.

Every chaotic region $R$ of $([\widetilde{f}_4];\widetilde{T}_4)$ is now a rectangle. 
Taking a minimal backwards stable iterate $([\widetilde{f}_5];\widetilde{T}_5)$ of $([\widetilde{f}_4];\widetilde{T}_4)$ gives a trellis mapping class such that every domain $D$ of $([\widetilde{f}_5];\widetilde{T}_5)$ with boundary in $\widetilde{T}_5^U\cup \widetilde{f}_5(\widetilde{T}_5^S)$ containing a chaotic region is a rectangle.
Foliate every domain $D$ containing a chaotic region, foliate $D$ and $f(D)$ by an unstable foliation $\cal{F}^U$ parallel to $T^U$ and a transverse stable foliation $\cal{F}^S$ parallel to $T^S$.
Isotope $\widetilde{f}_5$ to obtain a diffeomorphism $\widetilde{f}$ which preserves the stable and unstable foliations,
 and for which all points of non-chaotic regions are in the basin of a stable or unstable periodic orbit.
Let $(g;G,W)$ be a graph representative of $([\widetilde{f}_5];\widetilde{T}_5)$ for which $G$ is transverse to $\widetilde{T}_5^S$, and $\pi:\cutp{\widetilde{T}_5}\exto(G,W)$ be a deformation-retract which collapses each leaf of $\cal{F}^S$ onto a point of $G$.
Then $\pi\circ \widetilde{f}=g\circ\pi$ on every chaotic region $R$, so $\htop(\widetilde{f})=\htop(g)=\hniel[f;T]$.
Further, it is clear that every $\widetilde{f}$-extension of $\widetilde{T}_5$ is minimal, so every $\widetilde{f}$-extension of $T$ is minimal.
\end{proof}

%----------------------------------------------------------------------------------------------------------------------------------------------------------------------------

\subsection{Approximate entropy minimising diffeomorphisms}
\label{sec:minimisingsequence}

Throughout this section we take $([f];T)$ to be a well-formed irreducible trellis mapping class, $\lambda=\exp(\hniel[f;T])$ and $\lambda_\epsilon=\exp(\hniel[f;T]+\epsilon)$.
We aim to show that we can always find a diffeomorphism $\widehat{f}\in([f];T)$ which has entropy less than $\lambda_\epsilon$.
We first consider how to perform non-minimal extensions without increasing the Nielsen entropy above $\hniel[f;T]+\epsilon$.
To control the entropy bound, we need to consider the graph representative, which introduces some technical difficulties.
To avoid having to deal directly with graph representative in the sequel, we prove a result, Lemma~\ref{lem:transitiveextensionincrease}, which applies directly to trellis mapping classes.

Our final goal is to to construct a trellis mapping class which satisfies the conditions of Theorem~\ref{thm:minimalmodel}.
To do this, we may need to introduce new periodic points to the trellis to create attracting and repelling regions.
We then iterate curves bounding an attractor or repellor into regular domains, and finally take non-minimal iterates to move bigon boundaries into attractors and repellors.

%\fig{quadrantgraph}{A trellis and its graph representative in the neighbourhood of a periodic control edge $z_{p}$.}
%To show that we can take non-minimal iterates without increasing the Nielsen entropy, we consider a length function on a graph near a periodic control edge, as shown in \figref{quadrantgraph}.
%We let $q^u$ be the first intersection of $T^U(Q)$ with $T^S$, and $S$ be the stable segment on the $Q$-side of $T^U(Q)$ and $q^u$, as shown.
%Let $z_{p}$ be the control edge crossing $S(Q)$, and $z$ the control edge crossing $S$.
%We let $\beta(Q)$ be the curve from the periodic control vertex $v_p$ to the periodic control vertex $v_x$.
%Note that $\alpha$ is homotopic to the inclusion of the path $z_{p}\beta z_{x}$ in $M$.

\begin{lemma}
\label{lem:transitiveextensionincrease}
Let $([f];T)$ be a trellis mapping class and $Q$ a period-$n$ quadrant of $([f];T)$.
Let $q^u\in T^U(Q)$ be the endpoint of a stable segment $S$ on the same side of $T^U(Q)$ as $Q$,
 and let $q^s\in T^S(Q)$ be the endpoint of an unstable segment $U$ on the same side of $T^S(Q)$ as $Q$.
Then for any $\epsilon>0$ there is a minimal stable extension $([\widehat{f}];\widehat{T})$ such that $\hniel[\widehat{f};\widehat{T}]<\hniel[f;T]+\epsilon$,
 and an integer $k$ such that $\widehat{f}^{-kn}(S)\subset \widehat{T}^S$ and intersects $U$ at a point $\widehat{r}$
 such that $\{p,f^{-kn}(q^u),q^s,\widehat{r}\}$ are the vertices of a regular domain for $Q$.
\end{lemma}

\begin{proof}
Choose $\lambda$ and $\lambda_\epsilon$ such that $\hniel[f;T]<\lambda<\lambda_\epsilon\leq \hniel[f;T]+\epsilon$.
Let $(g;G,W)$ be the graph representative of $[f;T]$ and $l$ be a length function on $G$ such that for all edges $e$ of $G$, $l(g(e))<\lambda l(e)$.
Let $z_Q$ be the control edge crossing $S(Q)$, and $z_0$ be the control edge crossing $S$.
Let $\alpha$ be an edge-path starting at $z_Q$ and finishing at the first vertex $v$ between $z_Q$ and $z_0$, and let $\beta$ be the edge-path from $v$ to the end of $z_0$.
Since $l(g(\alpha))/\lambda_\epsilon < l(g(\alpha))/\lambda \leq l(\alpha)$,
 there exists $k$ such that $l(g(\alpha))/\lambda_\epsilon + 2l(\beta)/\lambda_\epsilon^{nk} < l(\alpha)$.

\fig{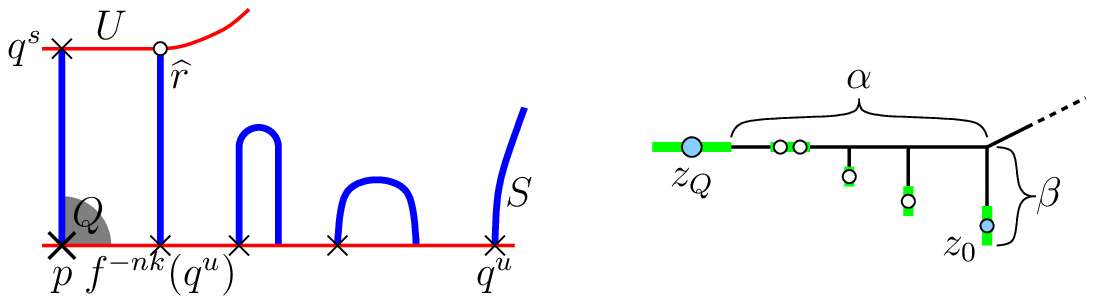}{A non-minimal extension.}
Let $\widetilde{f}$ be a diffeomorphism such that $\widetilde{f}^{1-kn}(T^S)$ is a minimal iterate of $T^S$, and set $\widetilde{T}=(T^U,\widetilde{f}^{1-kn}(T^S))$.
Then $[\widetilde{f};\widetilde{T}]$ has a graph representative $(\widetilde{g};\widetilde{G},\widetilde{W})$ for which $l$ extends to a length function with $l(\beta_{i})=l(\beta)/\lambda^{i}$ for $i$th backward iterates $\beta_i$ of $\beta$.
Let $\widehat{f}$ be a diffeomorphism which is isotopic to $\widetilde{f}$ relative to $\widetilde{T}$, for which $\widehat{f}^{-nk}(S)$ intersects $U$ transversely at a point $\widehat{r}$, but otherwise has minimal intersections, as shown in \figref{quadrantbifurcation}.
Let $\widehat{T}=(T^U,\widehat{f}^{-nk}(T^S))$, and $(\widehat{g};\widehat{G}\widehat{W})$ be the graph representative of $[\widehat{f};\widehat{T}]$.
The graph map $\widehat{g}$ maps $\alpha$ to $g(\alpha)$ and twice over $\beta_{nk-1}$,
 with total length $l(\widehat{g}(\alpha))=l(g(\alpha))+2l(\beta)/\lambda_{\epsilon}^{nk-1}$, so $l(\widehat{g}(\alpha)) < \lambda_\epsilon l(\alpha)$ as required.
Since the image of all other edges is unchanged, the growth rate of $l$ under $\widehat{g}$ is less than $\lambda_\epsilon$, so $\hniel[\widehat{f},\widehat{T}]=\htop(\widehat{g})<\log\lambda_\epsilon$.
\end{proof}

Before proving the main theorem we show how to ensure that every quadrant $Q$ is contained in a regular domain whose opposite sides bound attractors or repellors.
\begin{lemma}
Let $([f];T)$ be a transitive trellis mapping class, and $P$ be an essential periodic orbit of $([f];T)$ which does not shadow $T^P$.
Then there is a minimal supertrellis $([\widetilde{f}];\widetilde{T})$ of $([f];T)$ obtained by blowing up at $P$ such that every quadrant $Q$ is contained in a regular rectangular region such that $\widetilde{T}^S[q^u,r]$ is the boundary of a stable region $\widetilde{T}^S[q^s,r]$ is the boundary of an unstable region.
\end{lemma}

\begin{proof}
Since $P$ is an essential periodic orbit, every curve $\alpha(Q)$ eventually maps across $P$, so crosses a stable segment of $P$.
Hence iterating this segment backwards gives a curve bounding a regular domain, and further backward iterates give a regular region.
A similar argument holds for forward iterates of unstable segments, since in this case we can simply reverse time.
\end{proof}

\begin{theorem}[Existence of approximate entropy minimisers]
\label{thm:minimisingsequence}
Let $([f];T)$ be a well-formed trellis mapping class.
Then for every $\epsilon>0$, there exists a diffeomorphism $\widehat{f}\in([f];T)$ such that $\htop(\widehat{f})<\hniel[f;T]+\epsilon$.
\end{theorem}

\begin{proof}
We repeatedly construct supertrellises $([f_i];T_i)$ with $\hniel[f_i;T_i]<\hniel[f;T]+\epsilon=\log\lambda_{\epsilon}$.

Take $([f_0];T_0)=([f],T)$.
By irreducibility, there is a minimal extension $([f_1];T_1)$ of $([f_0];T_0)$ such that for every pair of nontrivial branches $T_1^U(p_u,b_u)$ and $T_1^S(p_s,b_s)$,
 there are points $p_u=p_0,p_1,\ldots,p_n=p_s$ such that $T_1^U(p_u,b_u)\cap T_1^S(p_1)\neq\emptyset$, $T_1^U(p_i)\cap T_1^S(p_{i+1})\neq\emptyset$ for $1\leq i<n-1$
 and $T_1^U(p_{n-1})\cap T_1^S(p_s,b_s)\neq\emptyset$.
By Lemma~\ref{lem:transitiveextensionincrease}, there is a (non-minimal) extension $([f_2];T_2)$ of $([f_1];T_1)$ such that every quadrant of $T_2$ is contained in a regular rectangular region.

We next construct a non-minimal extension $([f_3];T_3)$ of $([f_2];T_2)$ such that every nontrivial branch of $T_3^U$ intersects every nontrivial branch of $T_3^S$.
If $T^U(Q_U)$ intersects $T^S(Q)$ and $T^U(Q)$ intersects $T^S(Q_S)$, we can take minimal iterates of $T^U(Q_U)$ and $T^S(Q_S)$ until $T^U(Q_U)$ crosses $S(Q)$ and $T^S(Q_S)$ crosses $U(Q)$.
Then an application  of Lemma~\ref{lem:transitiveextensionincrease} shows that there is a non-minimal intersection such that $T^U(Q_U)\cap T^S(Q_S)\neq\emptyset$.

Now construct a non-minimal extension $([f_4];T_4)$ of $([f_3];T_3)$ such that there is an alpha-chain from any quadrant $Q_U$ of $T_4$ to any other quadrant $Q_S$.
Since $T_3^S(Q_S)$ intersects $T_3^U(Q_U)$, by taking minimal backward iterates of $T_3^S(Q_S)$ we can ensure that $T^S(Q_S)$ crosses $U(Q_U)$.
Then by Lemma~\ref{lem:transitiveextensionincrease}, there is a non-minimal extension such that $T^S(Q_S)$ crosses a regular region $R(Q_U)$.
Similarly, we can ensure $T^U(Q^U)$ crosses the regular region $R(Q_S)$.
A further application of Lemma~\ref{lem:transitiveextensionincrease} gives a non-minimal extension with an alpha-chain from $Q_U$ to $Q_S$.
Repeating this construction for all regions gives the required trellis $T_4$.

Let $([f_5];T_5)$ be a minimal supertrellis of $([f_4];T_4)$ such that all trivial branches of $T_4$ are contained in a non-chaotic region, as given by Theorem~\ref{thm:minimalsupertrellis}\ref{item:periodictrivialend} and \ref{item:trivialbranch}.
If $([f_5];T_5)$ has no attractors or no repellors, take a further minimal supertrellis $([f_6];T_6)$, with at least one attractor and one repellor, as given by Theorem~\ref{thm:minimalsupertrellis}\ref{item:periodicblowup}.
Since $([f_6];T_6)$ is transitive, we can take a minimal extension $([f_7];T_7)$ by iterating the stable and unstable curves bounding stable and unstable regions such that all quadrants of $T_7$ are contained in a regular region with the opposite stable side bounding an attracting domain, and the opposite unstable side bounding a repelling domain.

\fig{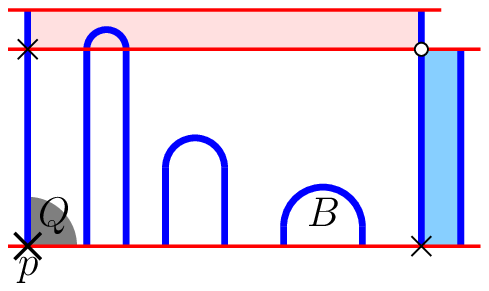}{Backward iterates of $B$ give an inner bigon in a repelling domain.}
Now, for any inner bigon $B$ such that $B^S$ is nonwandering, we take backward iterates of $B^S$ to obtain a stable segment $S$ in a regular domain of a quadrant $Q$.
Applying Lemma~\ref{lem:transitiveextensionincrease} gives a non-minimal extension such that $f^{-n}(B)$ contains a single inner bigon, and this inner bigon lies in a repelling domain, as shown in \figref{nonwanderingbigon}.
Similarly, by taking forward iterates of $B^U$ we can ensure $f^n(B)$ contains a single inner bigon in an attracting domain.
Applying this procedure to all bigons gives a trellis mapping class $([f_8];T_8)$.

Finally, we take a minimal supertrellis $([f_9];T_9)$ such that every chaotic region is a rectangle by introducing new stable curves as in Lemma~\ref{lem:invariantarcs}.
Then by Theorem~\ref{thm:minimalmodel} there is a uniformly-hyperbolic diffeomorphism $\widehat{f}\in([f_9];T_9)$, and hence in $([f];T)$ such that $\htop(\widehat{f})=\hniel[f_9;T_9]<\hniel[f;T]+\epsilon$ as required.
\end{proof}

We now prove the existence of pseudo-Anosov map which approximate the Nielsen entropy.
The condition on attractors and repellors is to ensure that the trellis mapping class contains a pseudo-Anosov map, since pseudo-Anosov maps have no attractors or repellors.
The strategy is to create the periodic orbits which will give the one-prong singularities of the pseudo-Anosov map.
Many of the steps of the proof mimic those the proof of Theorem~\ref{thm:minimisingsequence}.
\begin{theorem}[Existence of pseudo-Anosov representatives]
\label{thm:pseudoanosovmodel}
Let $([f];T)$ be a trellis mapping class with no attractors or repellors.
Then for any $\epsilon>0$ there exists a pseudo-Anosov diffeomorphism $\widehat{f}\in([f];T)$ such that $\htop(\widehat{f})<\htop([f];T)+\epsilon$.
\end{theorem}

\begin{proof}
If $T$ has an adjacent pair of trivial branches, take an extension $([f_1];T_1)$ such that these two branches intersect in a single transverse homoclinic point.
This does not affect the Nielsen entropy.
If $T$ has any other trivial branches, take a minimal extension $([f_2];T_2)$ such that these branches have an intersection point, which is possible by Theorem~\ref{thm:minimalsupertrellis}\ref{item:trivialbranch} since $T$ has no attracting or repelling regions.
As in the proof of Theorem~\ref{thm:minimisingsequence}, take an extension $([f_3];T_3)$ which is transitive.
Take a further extension $([f_4];T_4)$ such that every branch intersects every other with both orientations, which is possible by transitivity since we can take minimal iterates such that every unstable branch has an extension with and iterate of every stable branch which does not map to an intersection of $T_3$, and an isotopy in the neighbourhood of this intersection yields three intersections, one with the opposite orientation.

\fig{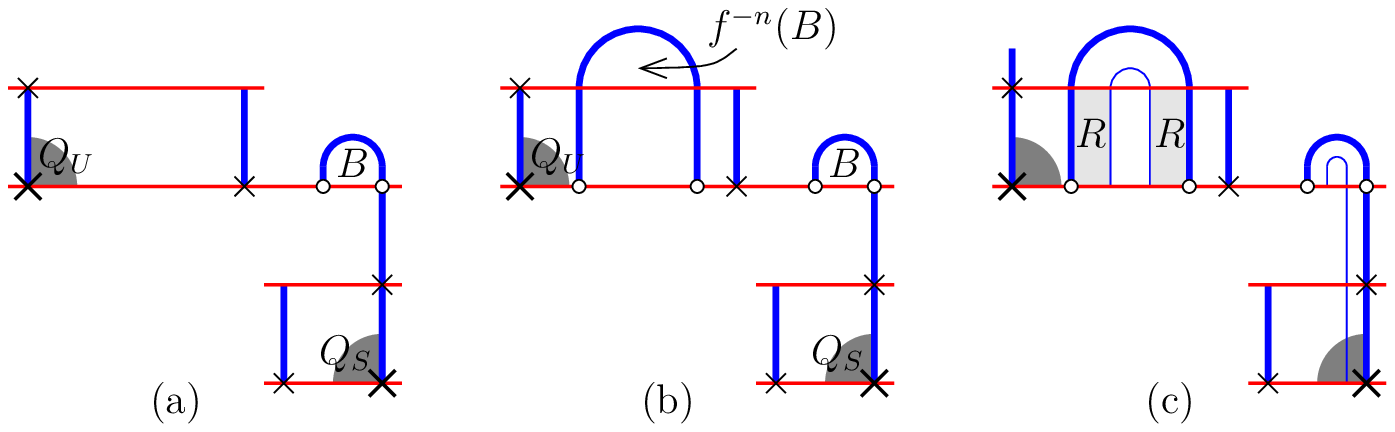}{The domain with vertices $q_0$ and $q_1$ must have periodic orbits in the rectangles $R_0$ and $R_1$.}
Let $B$ be an inner bigon of $T_4$, as shown in \figrefpart{chaoticbigon}{a}.
By Lemma~\ref{lem:transitiveextensionincrease}, we can find a stable extension such that $f{-n}(B^S)$ crosses a regular domain $R(Q_U)$ for some $n$.
Further, by removing intersections if necessary, we can ensure that $f{-n}(B^S)$ crosses $D$ twice and gives a new inner bigon $\widehat{B}$, as shown in \figrefpart{chaoticbigon}{b}.
A further application of Lemma~\ref{lem:transitiveextensionincrease} shows that we can find a stable extension $[\widehat{f}_3;\widehat{T}_3]$ such that $\widehat{B}^U$ crosses some regular domain, as shown in \figrefpart{chaoticbigon}{c}.
Since the regions $R$ are mapped over by $\alpha(Q_U)$ and map over $\alpha(Q_S)$, they must contain a periodic orbit, since $T$ is transitive.
Applying this construction for every inner bigon of $([f_4];T_4)$ gives a non-minimal extension $([f_5];T_5)$ such that every bigon of $T_4$ is a domain of $T_5$ containing an essential periodic orbit of $([f_5];T_5)$.

Let $([f_6];T_6)$ be the trellis mapping class obtained by puncturing at a periodic orbit in every inner bigon of $T^4$ to give a surface $M_6$.
Then $\hniel[f_6;T_6]=\hniel[f_5;T_5]$.
Since $T_4$ is a subtrellis of $T_6$, we can take a trellis mapping class $([f_7];T_7)=([f_6];T_4)$ in the surface $M_6$.
Since every inner bigon of $T_4$ contains component of $\partial M_6$, the trellis mapping class $([f_7];T_7)$ has no inner bigons.
The graph representative $(g_7;G_7,W_7)$ is locally injective except at cusps, so is efficient and hence is a train-track map for a pseudo-Anosov diffeomorphism $\widetilde{f}\in([f_7];T_7)$.
Then $\htop(\widetilde{f})=\hniel[f_7;T_7]\leq\hniel[f_6;T_6]<\hniel[f;T]+\epsilon$ as required.
\end{proof}

%----------------------------------------------------------------------------------------------------------------------------------------------------------------------------

\subsection{Non-existence of entropy minimisers}
\label{sec:noentropyminimiser}

Suppose that an irreducible trellis mapping class has a uniformly hyperbolic diffeomorphism realising the Nielsen entropy.
The following lemma shows that any isotopy removing intersections results in a trellis mapping class with strictly smaller Nielsen entropy.
If there is an alpha chain from any segment to any other, the trellis is said to be \emph{transitive}.

A transitive trellis has a transitive graph representative.
\begin{lemma}
\label{lem:transitivegraph}
Let $([f];T)$ be an transitive trellis mapping class.
Then the graph representative $(g;G,W)$ of $([f];T)$ has a single transitive component with positive topological entropy.
\end{lemma}
\begin{proof}
Let $\overline{G}=\bigcup_{n=0}^{\infty}g^n(G)$ be the essential subgraph, and let $e$ be any expanding edge of $G$.
By irreducibility, there is a quadrant $Q$ such that $g^n(e)$ tightens to $\alpha_{Q}$ for some $n$.

We now consider preimages of edges.
Again, let $e$ be an expanding edge of $G$.
There exists a tight curve $\beta_0$ in $\widetilde{G}$ with endpoints in $W$ such that $\beta_0$ is homotopic to a subinterval of a branch $T^U(p,b)$.
We then find a curve $\beta_1$ in $\widetilde{G}$ homotopic to a subinterval of a branch at $f^{-1}(p)$.
Proceeding recursively gives an edge-path $\beta_n=\alpha_Q$ where $\beta_i$ is a sub-path of $g(\beta_{i+1})$ for $0\leq i<n$.
Therefore, there exists a quadrant $Q$ and an integer $n$ such that $g^n(\alpha_Q)\supset e$. 

Now if $Q_U$ and $Q_S$ are any two quadrants, there exists $n$ such that $g^n(\alpha_{Q_U})$ contains a sub-curve $\alpha_{Q_S}$.
Since $(f;T)$ is irreducible, we can find $N$ such that for any two quadrants and any $n\geq N$, the iterate $g^n(\alpha_{Q_U})$ contains a sub-path $\alpha_{Q_S}$.
Therefore, any such curve $\alpha$ generates the same graph component under iteration.

Combining these results shows that there exists $N$ such that if $e_1$ and $e_2$ are any two edges, and $n\geq N$, then $e_2\subset g^n(e_1)$.
\end{proof}

\begin{theorem}
\label{thm:entropydecrease}
Let $([f];T)$ be a well-formed irreducible trellis mapping class, and suppose $f$ be a uniformly-hyperbolic diffeomorphism such that $\htop(f)=\hniel[f;T]$.
Then if $([\overline{f}];\overline{T})$ is a trellis mapping class which does not force $([f];T)$, then $\hniel[\overline{f};\overline{T}]<\hniel[f;T]$.
\end{theorem}
\begin{proof}
\fig{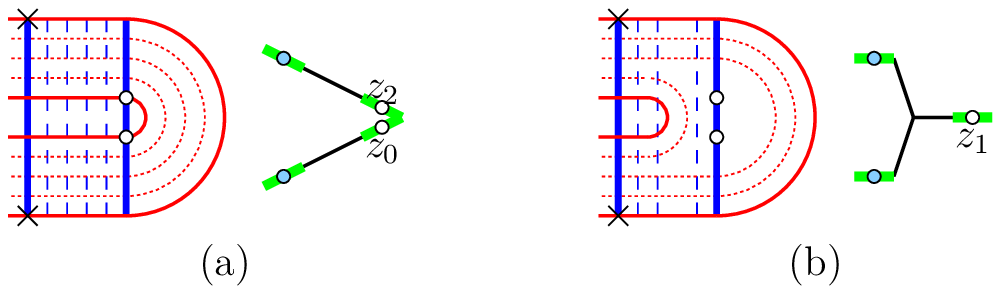}{Change in trellis and graph at a homoclinic bifurcation.}
By taking an $f$-extension $T_1$ of $T$, we can ensure that every inner bigon of $([f];T_1)$ is contained in a larger domain with the topology of \figrefpart{bifurcationdecrease}{a}, and that these domains are separate.
Further, we can ensure that every quadrant is contained in a rectangular region.
Now consider the effect of a homotopy to a trellis mapping class $([f_2];T_2)$ removing one pair of orbits on the same inner bigons.
The resulting trellis and graph locally have the topology of the right of \figrefpart{bifurcationdecrease}{b}.
The edges $a_0$ and $a_2$ must be expanding, since they forward iterate to a closed segment, and since every quadrant is contained in a rectangular region, $\partial g^n(a_0)=\partial g^n(a_2)=b$ for some edge $b$.
Therefore $a_0$ and $a_2$ are folded together to obtain the graph representative for $([f_2];T_2)$.
Further, there must be an edge mapping $\ldots\bar{a}_0a_2\ldots$.
Hence the entropy of the graph representative of $([f_2];T_2)$ is less than that of $[f;T]$, so $\hniel[f_2;T_2]<\hniel[f;T]$.
Since we can then prune $([f_2];T_2)$ to obtain $([\overline{f}];\overline{T})$, we have $\hniel[\overline{f};\overline{T}]\leq\hniel[f_2;T_2]$, hence $\hniel[\overline{f};\overline{T}]<\hniel[f;T]$ as required.
\end{proof}

We use this to prove that the sufficient condition given for the existence of entropy minimisers in Theorem~\ref{thm:minimalmodel} is necessary.
\begin{theorem}[Non-existence of entropy minimisers]
Let $([f];T)$ be an irreducible trellis mapping class.
Suppose that for every diffeomorphism $\widehat{f}\in([f];T)$ there is a $\widehat{f}$-extension of $T$ which is not minimal.
Then there does not exist a diffeomorphism in $([f];T)$ whose topological entropy equals $\hniel[f;T]$.
\end{theorem}

\begin{proof}
Suppose there is an $f$-extension $T_1$ of $T$ which is not minimal.
Let $T_2$ be a transitive $f$-extension of $T_1$, and $T_3$ be a further extension such that every non-wandering segment of $T$ crosses a regular domain of $T_3$.
Then we have entropy bound $\htop(f)\geq\hniel[f;T_3]$.

Take $([f_4];T_4)$ to be a minimal supertrellis of $([f];T_3)$ for which the opposite sides of every regular region bound a stable or unstable region.
Then every nonwandering segment of $T$ enters a stable or unstable region of $T_4$.
Prune $([f_4];T_4)$ to obtain a trellis mapping class $([f_5];T_5)$ which forces $([f];T_1)$ satisfying the conditions of Theorem~\ref{thm:minimalmodel}.
Then $\hniel[f_3;T_3]=\hniel[f_4;T_4]\geq\hniel[f_5;T_5]$, but by Lemma~\ref{thm:entropydecrease}, $\hniel([f_5];T_5)>\hniel[f;T]$.
Combining these inequalities we have $\htop(f)\geq\hniel[f_5;T_5]>\hniel[f;T]$ as required.
\end{proof}

The results of this section show that many fundamental properties of an irreducible trellis type $[f;T]$ depends on whether there is a diffeomorphism $\widetilde{f}\in[f;T]$ for which $\htop(\widetilde{f})=\hniel[f;T]$.
If such a diffeomorphism exists, then the entropy of the trellis type is carried in a uniformly hyperbolic diffeomorphism, but is fragile in the sense that any pruning will reduce the Nielsen entropy, and any diffeomorphism $\widehat{f}$ for which some extension is non-minimal must have strictly greater topological entropy.
If no such diffeomorphism exists, every diffeomorphism in the trellis type has an extension which is non-minimal, but pruning this extension gives a trellis type for which the entropy is still greater than the Nielsen entropy of $[f;T]$.

%****************************************************************************************************************************************************************************

%****************************************************************************************************************************************************************************

\bibliographystyle{alpha}
\bibliography{\bibdir/journalabbreviations,\bibdir/bibliography}

%****************************************************************************************************************************************************************************

\end{document}

%****************************************************************************************************************************************************************************

\section{Conclusions}
\label{sec:conclusions}

We have shown that the Nielsen entropy of a trellis mapping class is the infemum of the topological entropies in the class, and given necessary and sufficient conditions for this infemum to be a minimum.